\documentclass[11pt]{extarticle}

\pdfoutput=1

\usepackage{graphicx}

\usepackage{amsmath}
\usepackage{ragged2e}
\usepackage{amssymb} 
\usepackage{stmaryrd}
\usepackage{gensymb}

\usepackage{breakcites}

\usepackage{authblk}

\usepackage{geometry}
\usepackage{fancyhdr}
\makeatletter
\def\@biblabel#1{}
\makeatother

\newcommand{\citep}{\cite}

\title{Synergistic Dynamic Theory of Complex Coevolutionary Systems: Disentangling Nonlinear Spatiotemporal Controls on Precipitation}

\author[1]{Rui A. P. Perdig\~ao\footnote{Address: TU Wien, Karlsplatz 13/222, A-1040 Vienna, Austria. E-Mail: rui.perdigao@tuwien.ac.at.}}
\author[2]{Carlos A. L. Pires}
\author[1]{Julia Hall}
\affil[1]{Technische Universit\"at Wien (TU Wien), Austria }
\affil[2]{Instituto Dom Luiz, Universidade de Lisboa, Portugal}

\date{}

\begin{document}

\maketitle

\begin{abstract}

We formulate a nonlinear synergistic theory of coevolutionary systems, disentangling and explaining dynamic complexity in terms of fundamental processes for optimised data analysis and dynamic model design: Dynamic Source Analysis (DSA). 
DSA provides a nonlinear dynamical basis for spatiotemporal datasets or dynamical models, eliminating redundancies and expressing the system in terms of the smallest number of fundamental processes and interactions without loss of information.  This optimises model design in dynamical systems, expressing complex coevolution in simple synergistic terms, yielding physically meaningful spatial and temporal structures. These are extracted by spatiotemporal decomposition of nonlinearly interacting subspaces via the novel concept of a Spatiotemporal Coevolution Manifold. Physical consistency is ensured and mathematical ambiguities are avoided with fundamental principles on energy minimisation and entropy production. The relevance of DSA is illustrated by retrieving a non-redundant, synergistic set of nonlinear geophysical processes exerting control over precipitation in space and time over the Euro-Atlantic region. For that purpose, a nonlinear spatiotemporal basis is extracted from geopotential data fields, yielding two independent dynamic sources dominated respectively by meridional and zonal circulation gradients. These sources are decomposed into spatial and temporal structures corresponding to multiscale climate dynamics. The added value of nonlinear predictability is brought out in the geospatial evaluation and dynamic simulation of evolving precipitation distributions from the geophysical controls, using DSA-driven model building and implementation. The simulated precipitation is found to be in agreement with observational datasets, which they not only describe but also dynamically link and attribute in synergistic terms of the retrieved dynamic sources.

\end{abstract}

\section{Introduction}

\subsection{Fundamental Motivation}

The Earth System in general, and the Climate system in particular, are inherently coevolutionary and complex \citep{Peixoto+Oort1992, NicolisandNicolis2007}. Therefore, their fundamental understanding calls for the development and implementation of methodologies that disentangle nonlinear processes and interactions at play, ultimately improving their predictability. These are the core motivations of this paper. \\

The concepts of coevolution and complexity span a diversity of views depending on the background and context in which they are raised, including in biology, geomorphology and hydro-climatology, e.g. \cite{EhrlichandRaven1964, Whippleetal1999, PerdigaoBloeschl2014, Trochetal2015}. From all different accounts and perspectives, a general definition can now be formulated with overarching generality, with broad interdisciplinary applications in mind. 

A \textit{coevolutionary system} consists of a multiscale dynamical system wherein the intervening processes are dynamically bound to evolve in an interactive manner across spatiotemporal scales. Whether such systems are \textit{complex} depend on whether the system as a whole evolves in a manner that cannot be attributed to any individual or linearly combined factor, rather exhibiting features not present or explainable by any such factors alone: \textit{emergence}. 

Complexity thus stems neither from the variables themselves, nor from their linear combinations, but rather from their nonlinear interactions. If the intervening variables are independent from each other (i.e. do not have any dynamical mutual influence) and still interact to produce an emerging process, we are in the presence of a \textit{synergy}. For instance, secondary waves in triadic wave resonance are synergic processes emerging from nonlinear interacting pairwise independent waves \citep{PiresPerdigao2015} or oscillatory processes \citep{HockeandKaempfer2008}.

A wise statement from the ancient Greek philosopher Aristotle in his \textit{Metaphysica} compendia and highlighted in \cite{NicolisandNicolis2007} aptly makes the point on complexity: \textit{"The whole is more than the sum of its parts"}. 

The essence of the general methodological problem resides in finding the simplest set of non-redundant and dynamically independent controls able to explain and predict the dynamics of a physical process of interest in a general coevolutionary system. Given a set of measurements or signals capturing its behaviour, it is of interest to find and quantitatively characterise dominant features that represent the fundamental behaviour of the system, as represented by statistical or dynamical properties and interactions. This is particularly relevant in high-dimensional cases such as those involving atmospheric variables, where the full description of the system is often beyond reach and the dynamics appear to be random.
Albeit the apparent randomness, self-organised and synergistically sustained phenomena may emerge far from thermodynamic equilibrium \citep{Haken1983}, leading to coherent group behaviour that can be described by a simple set of macro-scale variables.  

\begin{figure}
\noindent\includegraphics[width=41pc]{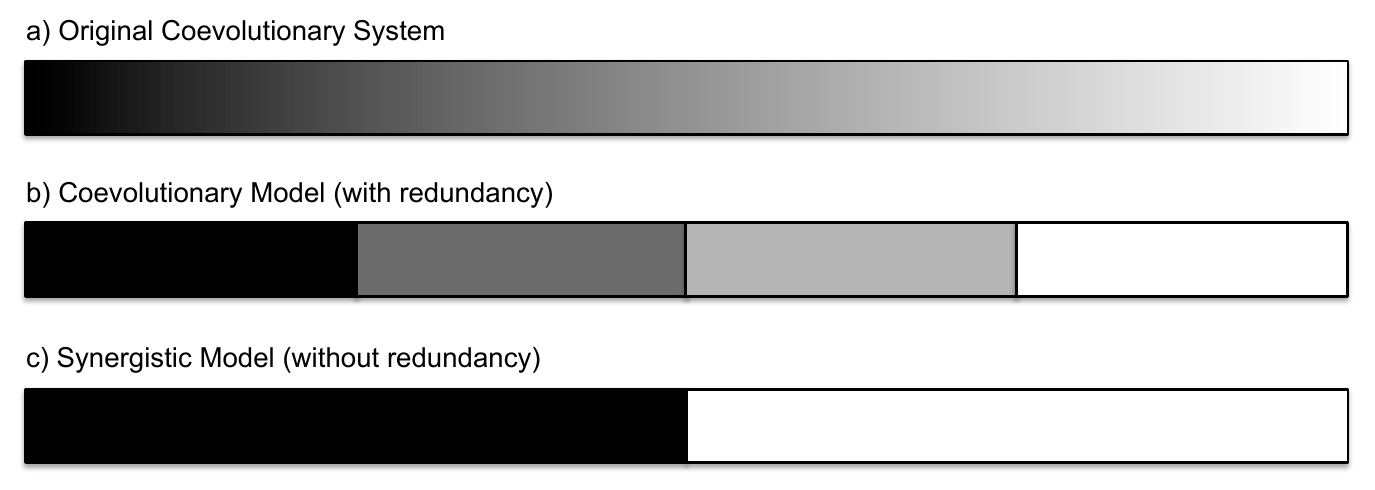}
\vspace{-10pt}
\caption{Schematic depiction of (a) a coevolutionary system, (b) a low-order coevolutionary model with some redundancy, and (c) the non-redundant synergistic model that we seek.}
\vspace{10pt}
\label{F1}
\end{figure}

Figure \ref{F1} depicts the nature of the problem in visually intuitive terms: (a) a complex coevolutionary system is represented by a broadband spectrum of gray shades with a large degree of redundancy as expressed by shading codependence; (b) a lower-order palette enables the reproduction of the original system using a few key components, whilst retaining some redundancy among shades; (c) the simplest set of non-redundant components consists of just two tones, black and white, the fundamental basis whose members can synergistically interact to generate the entire gray scale representing the original system. This is, in a nutshell, the fundamental set of basis components that we seek.

Bearing in mind the notion of a vector or functional space representing all possible states of the system, our goal is to find a basis for that space, i.e. the fundamental independent "axes" relative to which the entire dynamical process unfolds. As a functional basis, its components are independent from each other and their synergistic combinations generate the entire state space, i.e. the problem is reduced to the minimum possible dimensions without information loss.

In linear algebra and functional analysis, vector and functional spaces are traditionally generated by linear expansion of their basis elements \citep{ReedSimon1980, Roman2005}. Our interest, however, is in the more general case where the functional space is a spatiotemporal field fully generated by nonlinear combinations of an even lower-dimensional basis without loss of generality. For that purpose, we seek basis elements that are dynamically independent not only in the linear but also in the nonlinear sense. 
From within the multiplicity of mathematical solutions to this nonlinear problem, we further seek those actually representing fundamental physical processes. For that purpose, the sought basis terms are those for which first principles hold, namely on energy and entropy dynamics. 

All in all, we seek the smallest set of physically consistent, non-redundant dynamic processes with full generative power enabled by nonlinear dynamic interactions (i.e. a nonlinear synergistic basis), providing complex systems with an optimal model design in terms of fundamental processes and interactions.

\subsection{Applied Motivation}
While a general methodological undertaking is our main focus, a second, no less important motivation drives our quest: a dynamical understanding of Precipitation and its predictability in terms of underlying fundamental controls.

Precipitation is one of the most widely studied phenomena playing a crucial role in hydro-climate dynamics. Still, its thorough understanding and prediction remains rather elusive due to its dynamical complexity, stemming from nonlinear interactions among intervening processes.

In Statistical Climatology and Hydrology, Precipitation is often treated as a random process, with each event being a particular realisation of a theoretically or empirically inspired statistical distribution, conditioned or not by signatures of assumed controls in the causal chain.
Precipitation datasets are then often simulated as stationary or ciclo-stationary stochastic processes, \citep{WilksWilby1999, Paschalisetal2013}, or as random variables conditioned to relevant oceanic-atmospheric controls by Bayesian probabilistic models \citep{KidsonThompson1998, BardossyPegram2011, Songetal2014}, including nonlinear statistical downscaling approaches \citep{PiresPerdigao2007, Perdigao2010}.

Still in the statistical context, Precipitation datasets can also be studied with complex process network approaches \citep{Dongesetal2009, RuddellandKumar2009, Boersetal2013}, wherein intervening processes are interpreted as network nodes, and their relationships as network branches. In practice, the branches are simply statistical codependence measures ranging from traditional linear or nonlinear correlations to information-theoretical diagnostics. Therefore, while containing statistical value, they do not necessarily correspond to interactions in the physical sense.

Albeit the inherent complexity of Precipitation, its stochastic treatment can also be complemented by deterministic approaches of time series analysis, wherein self-similarity and scaling properties are brought out in connection with geometric features of the dynamics. The Fractal-Multifractal method \citep{Puente1992}, while inspired in the stochastic setting of multiscaling cascades in turbulence \citep{Mandelbrot1974}, is a deterministic-geometric approach already used for Precipitation studies \citep{Obregonetal2002, Maskeyetal2015}. 

Even though data-based analytics are aptly descriptive and quite practical to implement, they may leave out relevant physical information that is crucial to a better understanding and prediction.

\subsection{Traditional and Information-Theoretical Feature Extraction Methods}

While statistical feature extraction has a long tradition in the geosciences, the fundamental dynamical analysis we seek has been more elusive. 

Among various statistical techniques with diverse levels of sophistication, one of the most popular approaches to the extraction of spatial and/or temporal patterns is Principal Component Analysis (PCA), which has been widely used to extract dominant modes of spatial and temporal variability in atmospheric and oceanic fields \citep{Horel1981, WallaceandGutzler1981, Karletal1982, BarnstonandLivezey1987, Preisendorferetal1988}.
Fundamentally, PCA searches for uncorrelated components under a given inner product, while maximising the amount of variance they explain. The signal is then given as a linear combination of those components.

However, when the processes of interest have non-Gaussian or non-normal distributions, PCA and similar methods will give out uncorrelated principal components (PCs) that are not truly statistically independent. In fact, those components may be further decomposed as a nonlinear combination of a smaller number of underlying independent processes still able to generate the full signal. This problem is common to all factor analysis approaches in which the signal decomposition is done into uncorrelated components. It is worth reiterating that null correlation does not imply independence.

Statistical tools do exist to search for a non-redundant set of statistically independent features within a non-normal signal. One of the most notable examples comes from the Information-Theoretical framework on Independent Component Analysis (ICA) \citep{Comon1994, Hyvarinenetal2001} and equivalent methods. ICA searches for statistically independent components in the general case when the data is non-normally distributed. It does so by transforming a multidimensional dataset into a linear combination of components with minimised statistical dependence or Mutual Information from each other, i.e., as close to statistical independence from each other as possible.

Statistical independence means that the value of any of the components gives no statistical information on the values of other components, and therefore the statistical predictability potential of one variable from another one is null. 
Still, there will always be a small residual codependence among the components \citep{CoverandThomas1991}, with the theoretical foundations for a strictly positive lower bound for Mutual Information associated to prescribed cross moments having been derived by \cite{PiresPerdigao2012} independently of data size constraints, and by \cite{PiresPerdigao2013} for finite-sized data samples.

Within the geophysical sciences, ICA has been applied to the Sea Surface Temperature (SST) field extracting statistically independent modes of variability in \cite{Airesetal1999} and an intuitive didactic approach to ICA applied to geophysical-like correlated data was provided in \cite{Airesetal2002}.  Moreover, the method was also used by \cite{FodorandKamath2003} to separate dominant statistical features in climate data, in \cite{Basaketal2004} for weather data mining using the North Atlantic Oscillation as a specific example, and by \cite{Perdigao2004} to produce statistically independent spatiotemporal teleconnection patterns from non-normally distributed atmospheric fields, using then those components as predictors for meteorological regimes.

For all its merits, ICA still has some limitations: its basic form does not disentangle nonlinearly mixed signals and its results are ambiguous in some situations, such as under rotational degeneracy when multi-dimensional subspaces are jointly Normally distributed, or when nonlinear mixing comes into play. 

These limitations may be addressed to some extent by considering the full signal as a mixing of multivariate non-normally distributed sources through the Independent Subspace Analysis \citep{PiresRibeiro2016}. Alternatively, nonlinear generalisations of ICA and PCA can be considered \citep{Pajunenetal1996, Oja1997}, along with advanced methods for nonlinear statistical decomposition of kinematic features such as Principal Interacting Patterns (PIP) and Dynamic Mode Decomposition (DMD) \cite{Hasselmann1988, Tuetal2014}. 

However, the sophistication of these approaches often comes at the expense of the uniqueness and fundamental meaning of the solutions. Moreover, whilst capturing nonlinear statistical dependence among system features on an aggregate level, they do not capture the dynamical dependence within the spatiotemporal domain of analysis, i.e. the small scale behaviour beneath the overall large scale features captured by the statistics. This is a common aspect of all statistical techniques involving \textit{blind source separation} \citep{Yuetal2014}.

The core of the problem is that \textit{processes that are statistically independent can still be dynamically codependent} \citep{Perdigao2016}. By leaving out the fine print in the dynamics, the aggregate diagnostics have a more limited predictive potential. Even if in the end an aggregation is done on dynamical diagnostics, at least their fundamental structure will have "learnt" (gathered information) from the dynamics rather than solely from the aggregate statistical features.

\subsection{Study Outline}

One of the key aims of the present study is to disentangle nonlinearly interacting processes in the climate system, in order to bring out fundamental mechanisms that, whilst dynamically independent from each other, cooperatively influence precipitation regimes.
For that purpose, this paper introduces a physical source decomposition methodology in section \ref{s2}, where the fundamentals of dynamic source analysis are formulated. 
A framework for nonlinear space-time decomposition in general spatiotemporally coevolving fields is then introduced in section \ref{s3}. The theoretical foundations for evaluation of predictability and dynamic model building are then introduced in section \ref{s4}. \textbf{The methodologies introduced in the paper are essentially nonlinear dynamic analysis and model building techniques with broad applications.}
 The application of the novel methodologies to the retrieval of nonlinear atmospheric controls of precipitation is conducted in section \ref{s5}, followed by an evaluation of dynamic predictability and simulation of evolving precipitation distributions given the retrieved controls, in section \ref{s6}. The main body of the paper is then completed with an overall discussion and concluding remarks in section \ref{s7}. 
 
Further to the main body of the paper, three appendices are included. Appendix \ref{AA} addresses the concepts of coevolution and synergy in information-theoretical terms. Appendices \ref{AB} and \ref{AC} introduce new theoretical developments in functional analysis necessary for the spatiotemporal decomposition conducted in the dynamic source analysis framework. These methodologies enable an effective decomposition of nonlinearly interacting functionals (e.g. the spatial and temporal information) even in a coevolutionary setting whereby space and time are nonlinearly entangled by relative celerities.

\section{Dynamic Source Analysis: a Synergistic Theory of Coevolutionary Systems}\label{s2} 

\subsection{Motivation and Fundamentals}\label{s21}

\subsubsection{Expressing Redundant Coevolution as Non-Redundant Synergy}\label{s211}

Coevolutionary complex systems involve codependent observables, i.e. mensurable quantities with some degree of redundancy. This is indeed the case in the current paradigms of coupled dynamical systems. The existence of redundant information suggests that a more efficient formulation could be derived involving a lower number of fundamental non-redundant terms without loss of generality. The fundamental terms are the independent dynamic sources we seek, representing the fundamental "backbone" of the dynamical system. 

While coevolution always entails a degree of redundancy among processes, a synergistic interaction does not. In fact, it rather entails cooperative dynamics wherein the intervening processes jointly produce emerging structures with features that do not exist in any individual or linearly combined source components.
These matters are mathematically exemplified in Appendix \ref{AA}. 

A purely synergistic interaction is optimally efficient since it leads to the emergence of unprecedented features (\textit{innovation}), giving the ability of a simple set of processes to synergistically produce an emergent cascade of child processes that ultimately leads to the full observed complexity \citep{PiresPerdigao2015}. The polyadic wave resonance in fluids is a simple example of synergistic interactions among parent processes ($i^\text{th}$ generation waves) leading to the emergence of child processes [$(i+1)^\text{th}$ generation waves]. This is the case in geophysical phenomena such as nonlinear triadic and quartic interactions among ocean waves \citep{Komenetal1996}.

Essentially, our idea is to express coevolutionary complex systems involving a large number of codependent observables in terms of synergistic interactions among independent components representing the fundamental underlying processes. 
For that purpose, we seek a simple nonlinear dynamic basis for complex nonlinear systems. 

\subsubsection{Unveiling a Simple Nonlinear Basis for Complex Systems}\label{s212}

Let $\mathbf Y_{s,t}$ denote the set of observables, i.e. the physical variables or fields that are observed, measured and recorded in spatiotemporal datasets, where $s$ represents the spatial coordinates and $t$ the time. Let $\mathbf X_{s,t}$ denote the set of the hereafter called \textit{Dynamic Sources}, taken as spatiotemporal functions from a Hilbert space \citep{CourantHilbert1953} where differentials can be defined: a Sobolev space \cite{Sobolev1938, Sobolev1963} $\mathcal H^\beta$ comprising a class of $\mathcal L^2$ (square integrable) functions with weak derivatives up to order $\beta$. Weak differentiability means that differential operators can be defined everywhere in the domain even if traditional differentiability fails in a null-measure set (e.g. a 2D set in a 3D domain). This is particularly useful when evaluating fluid flows over critical transitions, e.g. across weather fronts \citep{Perdigao2016}.

Essentially, the dynamic sources are fundamental spatiotemporal state variables from which $\mathbf Y_{s,t}$ depend via an operational set of dynamical relationships $\mathbf f$ in the general form:
\begin{equation}\label{e1}
\mathbf Y_{s,t} = \mathbf f(\mathbf X_{s,t})
\end{equation}

The operator $\bf f$ maps the $m$-dimensional functional domain of the sources, $D(\mathbf X_{s,t})$ to the $n$-dimensional functional domain of the observables, $D(\mathbf Y_{s,t})$, both of which are subsets of $\mathcal H^\beta$. In practice, $\bf f$ can be seen as an \textit{observation operator} akin to those involving spatial re-colocation, temporal delaying, differentiation or integration. Without a priori constraints to the spatial and temporal domains of sources and observations, they are not necessarily the same. Moreover, the spatial coordinates may be time-dependent in a coevolutionary setting, leading to a celerity term $ds/dt$ of the space referential \citep{PerdigaoBloeschl2014}.

The aforementioned variables are generic $m$-dimensional functions living in an $n$-dimensional functional space, with $m\leq n$. The dynamical system in equation \eqref{e1} then prescribes the dynamics of a spatiotemporally distributed quantity, which can be straightforwardly treated as a multidimensional deterministic field or a multivariate stochastic variable or distribution.

Our goal is to retrieve the dynamic sources $\mathbf X_{s,t}$ from the observational datasets $\mathbf Y_{s,t}$ that minimise the dynamic dependence (linear and non-linear) among the components of $\mathbf X_{s,t}$ and maximise their ability to [not necessarily linearly but also nonlinearly] generate the entire state space spun by $\mathbf Y_{s,t}$. In other words, given the observations (the observed dynamics) we seek a functional basis in their state space. That basis shall provide the optimal components underlying the observations, since it represents them with neither redundancies nor information loss.

The retrieval of $\mathbf X_{s,t}$ from $\mathbf Y_{s,t}$ would be trivial if their functional relationship $\bf f$ would be a known invertible function, by simply taking $\mathbf X_{s,t} = \mathbf f^{-1} (\mathbf Y_{s,t})$. However, that is not necessarily the case in most applications. 
Without further knowledge, the problem would be under-determinate. Therefore, constraints could be put into place, as done statistically with ICA and PCA, so that ambiguities would be overcome.  In the present study, the constraints are of physical nature (section \ref{s25}) and the feature extraction performed over the dynamics (sections \ref{s22} - \ref{s24}). For notational ease, unless otherwise specified, symbols in uppercase bold type shall refer to spatiotemporal variables, i.e. $\mathbf X \equiv \mathbf X_{s,t}$ and $\mathbf Y \equiv \mathbf Y_{s,t}$ in the upcoming sections.

\begin{figure}
\noindent\includegraphics[width=40.5pc]{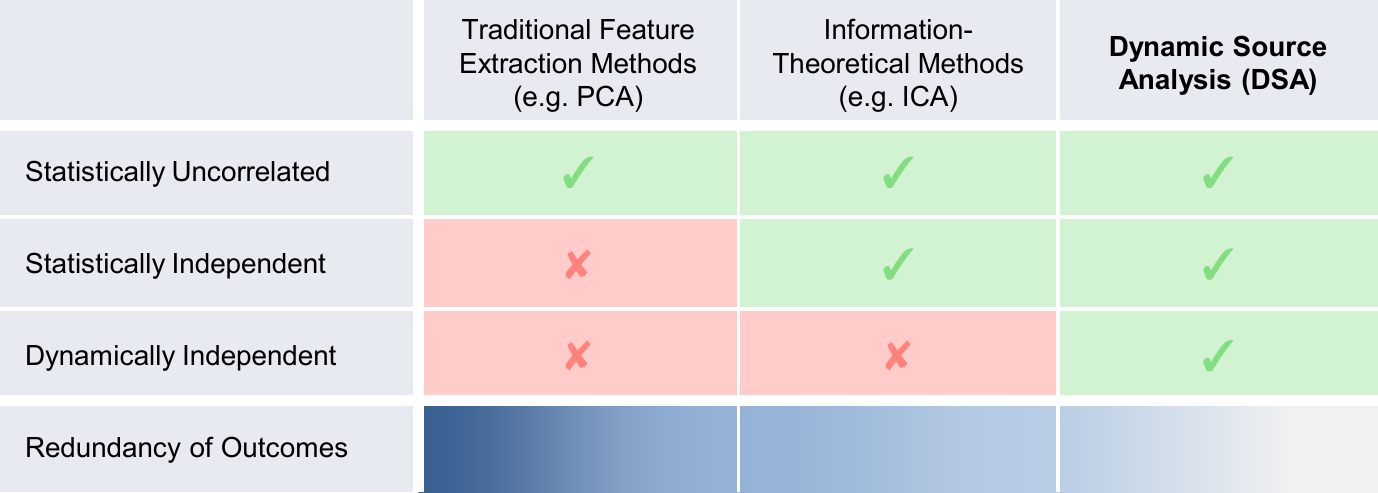}
\vspace{-10pt}
\caption{Comparison among Feature Extraction Methods in terms of the ability to extract independent components and thus minimise the number of factors needed to fully characterise the system.}
\vspace{10pt}
\label{F2}
\end{figure}

The key differences between the proposed Dynamic Source Analysis (DSA) and existing families of methods are summarised in Figure \ref{F2}. Essentially, traditional feature extraction methods yield uncorrelated but not statistically independent components; information-theoretical methods do provide statistically independent components but not dynamically independent ones; and DSA provides dynamically independent components, thus yielding the lowest level of redundancy in the characterisation of the system.

\break

\subsection{Nonlinear Diagonalisation of Dynamic Interactions}\label{s22}

Dynamic Source Analysis seeks dynamically independent sources, i.e. such that no redundant dynamic interactions exist among them. This brings out a canonic formulation of the dynamical system in terms of non-redundant synergistic interactions, which is naturally more efficient than formulating the dynamics in redundant, coevolutionary terms.

For this purpose, we define \textit{Dynamic Interaction of $k^{th}$ order} among the components of $\mathbf X$ as:
\begin{equation}\label{eq2}
\mathcal D^k(\mathbf X) \equiv
{\nabla^k_{\bf X} \, {\bf \dot X}}
\end{equation}
where $\nabla^k_{\bf X}$ is the $k^{th}$-order tensor differential operator given by the $k$-times recursive application of the gradient with respect to the spatial components of $\mathbf X$, and $\bf \dot X$ is the time derivative of $\mathbf X$.
In the spectral case where each spatiotemporal component $X_\lambda$ corresponds to a different scale $\lambda$, equation \eqref{eq2} provides the \textit{Dynamic Scale Interaction of $k^{th}$ order} within $\mathbf X$.

The dynamically independent sources are thus those for which the following condition holds:
\begin{equation}\label{dnull}
\mathcal D^k(\mathbf X) = \mathbf \Lambda_k, \,\,\,\, \forall \, {k}
\end{equation}
where $\mathbf \Lambda_k = \text{diag} (\Lambda_{k,1} ; \cdots ; \Lambda_{k,m})$, with $m=\text{dim}(\mathbf X)$, is a functional basis to the dynamical system prescribed by $\bf \dot X$. As such, the off-diagonal elements of $\mathcal D^k(\mathbf X)$ vanish, meaning that the cross-dependencies are null not only at bilateral $(k=2)$ but also multilateral $(k>2)$ levels. In the linear case, $\mathcal D^k = 0$ for $k>1$ and $\mathcal D^1$ corresponds to a diagonal Jacobian matrix of the underlying dynamical system, i.e. null linear couplings.
In the general nonlinear case, it means that along the whole trajectory in state space $\mathcal S$ all non-diagonal (i.e. cross) dynamical interactions are null.

The derivatives in equation \eqref{eq2} require the knowledge of the phase space trajectory prescribed by $\bf \dot X$. However, that information can be found from the observable dynamics $\bf \dot Y$. In fact, since both $\bf \dot Y$ and $\bf \dot X$ come from the same dynamical system, they have a topologically equivalent phase space portrait except if $\mathbf X$ in equation \eqref{e1} is under-determined from $\mathbf Y$. 

In practice, the dynamical systems of interest shall be generic spatiotemporal fields spanning phase space \textit{manifolds}. These are topological spaces that, whilst eventually curved, can be mapped, at the vicinity of each point, to and from Euclidean spaces via a functional bijection (invertible one-on-one correspondence or isomorphism). This is the case with fundamental geodesic and relativistic spatiotemporal geometries with broad applications.

One of the most relevant isomorphisms is the transformation between non-Euclidean manifolds and the corresponding Euclidean tangent spaces at the vicinity of each point, which under appropriate differentiability conditions defines a \textit{diffeomorphism} \citep{DaSilva2001}. Therefore, albeit the curvature of the spatiotemporal field of study, the topological space spun by its spatiotemporal gradients is conveniently comprised of an Euclidean structure (\textit{tangent bundle}, union of all local tangent spaces).

A crude motivational image of a tangent bundle can be obtained by depicting a ball covered by flat "postit" stickers. The ball represents the manifold, each sticker is a local tangent space, and the family of stickers covering the entire ball is the tangent bundle.

The rigorous mathematical treatment of manifolds and associated concepts stems from  \textit{differential geometry}, whose roots date centuries back to the origins of cartography \citep{Spivak2005}.

\subsection{Beyond Statistical Independence: Differential Geometric Quest for Dynamic Independence}\label{s23}

Essentially, DSA can be seen as a dynamic generalisation of ICA wherein the independent components are computed at the tangent space to each point and fundamental physical constraints on energy and entropy dynamics (section \ref{s25}) are taken into consideration. 
The procedure is conducted throughout the phase space manifold to yield the global dynamic sources. Given the differential geometric isomorphism between local charts and tangent spaces, the independent components computed on the tangent spaces are applicable to the actual manifold chart as well. The global dynamic sources are thus defined in the tangent bundle guiding the underlying phase space manifold. 

With this procedure, DSA ensures independence among components at every stage of the dynamics (dynamic independence), rather than only over a statistical aggregate (statistical independence as in ICA and information theory). 
Moreover, the dynamic information is preserved as the analysis follows the phase space trajectories rather than assigning classes with disregard to the dynamic sequence (as happens when defining a statistical distribution e.g. by assigning bins or classes).

Albeit the technical sophistication of differential geometry, it should be noted that its mathematical essence is actually rather intuitive if we think in geodesic terms. 
In all quantitative geosciences, local coordinates (e.g. zonal, meridional) are also defined on the local tangent space to the curved geophysical surface of interest at each point. While the local coordinate system evolves from point to point due to the curvature of the surface, a global set of coordinates is then defined with overarching span and simplicity, e.g. the spherical coordinates. Similarly, the dynamic sources are fundamentally global functions that encompass, in a similarly simple manner, the entire phase space manifold. As such, they define a basis grounded on a canonic set of generalised curvilinear coordinates. The practical challenge thus lies in finding such a functional basis for generic dynamical systems.

\subsection{Dynamically Independent Sources from Coevolving Observables}\label{s24}

The problem consists essentially on expressing a complex coevolutionary dynamical system $\bf \dot Y$ in Taylor expansion terms of a simpler synergistic system of independent dynamic sources $\bf \dot X$:
\begin{equation}\label{e3}
{\bf \dot Y} = \sum_{k} k!^{-1} \left[ {\nabla^k_{\bf \dot X}{\bf \dot Y}}  \right]_{\mathcal R} \, {\bf \dot X}^k
\end{equation}
where $\mathcal R$ refers to evaluation at the reference manifold\footnote{The reference manifold is spun by the invariants of motion and includes not only attractor but also repeller sub manifolds.}, and
\begin{equation}\label{e4}
{\bf \dot X} = \sum_k k!^{-1} \left[{\mathcal D^k (\mathbf X)}\right]_\mathcal{R} \, \mathbf{X}^k ,
\end{equation}
where $\mathbf X$ are the spatiotemporal fields characterising the dynamic sources ($\mathbf X \equiv \mathbf X_{s,t}$ as defined in Section \ref{s212}).

By applying the condition \eqref{dnull} to equation \eqref{e4}, the dynamically independent sources sought after are ultimately given by:
\begin{equation}\label{e4di}
{\bf \dot X} = \sum_k k!^{-1} [{\bf \Lambda}_k]_\mathcal{R} \, \mathbf{X}^k \, .
\end{equation}

In practice, a functional diagonalisation is performed to the coevolutionary system towards meeting the condition in \eqref{dnull} and thus leading to \eqref{e4di}.
For that purpose, DSA searches for $\mathbf X$ that minimises the off-diagonal functionals defined by

\begin{equation}\label{qn}
\nu_k (\mathbf X) = || \mathcal D^k(\mathbf X) -  \mathbf \Lambda_k ||^2 \, , \,\,\,\, \forall k
\end{equation}
with $|| \cdot ||^2 \equiv \langle \cdot, \cdot \rangle$, 
i.e. the quadratic norm or self inner product of the difference between the left- and right-hand sides of the dynamical independence condition \eqref{dnull}. 

The quadratic norm in equation \eqref{qn} is taken in the Sobolev space $\mathcal H^\beta$ with the inner product generically defined by:
\begin{align}\label{innerprod}
\langle g, h \rangle_{\beta,\Theta} 
&= \sum_{\alpha \leqslant \beta} \langle \mathrm{d}^\alpha g, \mathrm{d}^\alpha h \rangle_{\mathcal L^2(\Theta)} \\
&= \sum_{\alpha \leqslant \beta} \int_\Theta [\mathrm{d}^\alpha g]^* \, \mathrm{d}^\alpha h \, d \mu
\end{align}
where $g, h \in \mathcal L^2(\Theta)$ are generic square-integrable tensor functions, $\beta$ is the order up to which differentials are defined in the Sobolev space, $\alpha \in\{0,\cdots,\beta \}$ is the order of the differentiation operator $\mathrm d$, $\Theta$ is an
unbounded manifold living in $\mathcal H^\beta$, and * denotes \textit{complex conjugate transpose}. The measure $\mu$ refers in general to any mathematically consistent metric in $\mathcal H^\beta$.

Note that, while coevolutionary dynamics unfold among the components of $\mathbf Y$, no such coevolution occurs among different independent dynamic sources $X_i$ in $\mathbf X$. In fact, by construction the components $X_i$ interact in a cooperative synergistic manner in the dynamics of $\mathbf Y$ whilst retaining their own dynamic independence relative to each other.

\subsection{Disambiguation by Physical Principles}\label{s25}

Eventual ambiguities in the solutions stemming from non-linearities (e.g. solutions of $\mathbf X$ with the same quadratic norm) are overcome by taking the physical solution grounded on first principles, namely the dynamic configuration with maximum entropy production and minimum energy consumption. Entropy and energy are quantified in the mechanistic sense, which enables their characterisation in dynamical systems terms by resorting to topologic properties of the phase space representing the dynamics. In practice, the entropy measure is computed with the \textit{topologic entropy} $h_\text{top}(\mathbf X)$ \citep{Ott2002} and the energy as a \textit{Morse-Lyapunov function} $\phi(\mathbf X)$ \citep{GrinesPochinka2010}. \\

We caution the reader that purely dynamical systems approaches are only kinematic or motion-descriptive, until adequately taking physical principles into consideration. It is important to take into consideration that the fundamental entropy production principles hold at constant energy levels, and energy minimisation principles at constant entropies \citep{Lage1995}. Therefore, our fundamental physical constraints consist of seeking solutions with \textit{minimum isentropic energy} and \textit{maximum equienergetic entropy production}.  That is, entropy production is evaluated in the equienergetic subspace and energy dynamics evaluated in the isentropic subspace of the phase space spun by the dynamical system. The joint implementation of these principles brings out the general disambiguated sources over the overall phase space. Moreover, it ensures the physical consistency of the solutions. 

Formally, we define the \textit{Physically Optimised System} $\Psi(.)$ among the functional of mathematical possibilities $\Omega$ for equienergetic entropy production rates $\gamma$ and isentropic energy consumption rates $\xi$, as the concrete function within $\Omega$ for which the equienergetic entropy production rate is maximum and isentropic energy consumption rates is minimum: 

\begin{equation}
\Psi(\Omega) \equiv \delta[\Omega_{(\gamma,\xi)},(\gamma_\text{max},\xi_\text{min})]
\end{equation}
where the Delta functional $\delta[\mathbf f(\mathbf x),\mathbf a]$ extracts the concrete characterisation of the function $\bf f(x)$ at $\bf x=a$.

By applying this definition to the set of mathematical solutions $\Omega(\mathbf X_\text{math})$, the physically optimised system $\mathbf X$ becomes:
\begin{equation}
\mathbf X = \Psi[\Omega(\mathbf X_\text{math})] = \delta[\Omega(\mathbf X_\text{math})_{(\gamma,\xi)},(\gamma_\text{max},\xi_\text{min})]
\end{equation}

Having retrieved dynamically independent sources from the nonlinearly mixed observables captured by the measurements, it is important to ensure the invariance of the information content within the system, i.e. that no information has been lost or spuriously added in the procedure.

In the dynamic source analysis, the dynamics of the measurements $\bf Y$ and of their corresponding independent sources $\bf X$ are related by a locally smooth homeomorphism, equation \eqref{e3}. The physical reason for this is that their respective dynamical systems, $\bf \dot Y$ and $\bf \dot X$, span the same phase space manifold and invariants of motion, as noted above. Therefore, they also share the same information content, which thus ensures the non-existence of spurious or lost information in the analysis. That is, the referential may change but the underlying physics do not.

\section{Nonlinear Space-Time Decomposition and the Spatiotemporal Coevolution Manifold}\label{s3}

\subsection{Motus: Expressing Spatiotemporal Processes in Observable Spatial and Temporal Structures}\label{s31}

Having retrieved dynamically independent processes, it is also important to contextualise, characterise and depict them in terms of perceptual references such as space and time. 
In that regard, we will be interested in their spatial structure captured over some time span (e.g. a \textit{Spatial Structure over a Climatology}); and on the temporal structure captured over some spatial domain (e.g. a \textit{Temporal Structure over a Region)}. 
For that purpose, observable spatial and temporal structures are sought for the dynamic sources $\mathbf X$.

A standard procedure in space-time decomposition consists of taking subspace projections of the spatiotemporal process by statistically aggregating it over the orthogonal complement of the subspace of interest. For instance, in a linearly separable space-time coordinate system, spatial patterns can be obtained by averaging the spatiotemporal process over time, and aggregate time series by averaging the process over space. A spatiotemporal process can then be expressed in terms of its temporally and spatially distributed components, e.g. Principal Components (time) and Empirical Orthogonal Functions (space) in Principal Component Analysis. Given their statistical nature, they assist in the study of its temporal and spatial variability.

However, in the aforementioned procedure the dynamics are lost due to the averaging operations. Moreover, spatial and temporal subspaces are not necessarily independent. In fact, the spatial reference can be time-dependent, as expressed via spatiotemporal coevolution \citep{PerdigaoBloeschl2014}.

\subsection{Decomposition over Interacting Subdomains}

In this section, a dynamically preserving methodology is introduced for space-time decomposition of spatiotemporal processes in the general case where spatial and temporal components of the reference frame are not necessarily independent from each other. 

Consider a generic functional basis $\mathbf X_{s,t}$ which can be the independent dynamic sources retrieved in the previous section. Whilst independent from each other, the components $(X_i)_{s,t}$ of $\mathbf X_{s,t}$ are all spatiotemporal dynamical processes. Our aim is to characterise these spatiotemporal components into observable spatial and temporal functionals. 

For that purpose, we introduce an \textit{Interacting Subspace Decomposition Operator} $\Upsilon_l(.)$ that, once applied to $\mathbf X_{s,t}$, extracts its spatial and temporal subspace manifolds (sub manifolds) for $l=s, l=t$ respectively. 
The general mathematical definition of the operator $\Upsilon_l(.)$ is introduced in Appendix \ref{AB}, and its application to $\mathbf X_{s,t}$ is given by:
\begin{equation}\label{op}
\Upsilon_l(\mathbf X_{s,t}) \equiv \mathbf X_{s,t} \star \mathbf e_l  \equiv \mathbf x_l
\end{equation}
where we introduce a \textit{Retrieval Product} $A \star B$ retrieving the $B$-structure or dimension from $A$. Its general definition is also introduced in the Appendix \ref{AB}. Still in equation \eqref{op}, $\mathbf e_l$ is the normalised functional basis for the subspace $l$, and $\mathbf x_l$ is a tensor functional with the rank $r_l$ of the projected subspace.
For instance, $r_s=3$ and $r_t=1$ for a tridimensional space and unidimensional time.

The spatiotemporal decomposition in \eqref{op} is performed whilst concentrating all the space-time codependencies into a generalised coevolution term, consisting of the interactions between spatial and temporal submanifolds. In order to mathematically characterise this term, we introduce the concept of a \textit{Spatiotemporal Coevolution Manifold}, hereby defined as:
\begin{equation}\label{coevmanif}
\mathbf C_{s,t} \equiv \langle \mathbf e_s , \mathbf e_t \rangle
\end{equation}
where $\langle \cdot , \cdot \rangle$ represents the inner product defined in equation \eqref{innerprod}, involving multiple-order interactions. 
By noting that the gradient of a generic coordinate $l$ is given by $\nabla_{s,t}(l) = \mathbf e_l = \frac{\partial l}{\partial s} \mathbf e_s + \frac{\partial l}{\partial t} \mathbf e_t$, the first-order gradient differential form of the spatiotemporal coevolution is given in terms of relative celerities: $\mathbf C_{s,t}^{[1]} = \frac{\partial t}{\partial s} \mathbf e_s + \frac{\partial s}{\partial t} \mathbf e_t$. By definition of gradient, $\mathbf C_{s,t}^{[1]}$ is locally orthogonal to $\mathbf C_{s,t}$.

In spatiotemporally separable systems, $dim(\mathbf C_{s,t}) = 0$ as there is no shared dimension between space and time. In non-separable systems (in which spatial patterns are time dependent), the tensor rank of $\mathbf C_{s,t}$ characterises the dimensionality of the spatiotemporal coevolution manifold.

A simple example is now given for a propagating wave $W_{s,t}=\cos(s-v\,t)$ on a coevolutionary spatiotemporal coordinate system $s=\omega \,t$, with $v$ and $\omega$ as constants. Here, the coevolution manifold is the 1D manifold given by the straight line defined by $s-\omega \, t = 0$ and the magnitude of the spatiotemporal interaction is given by the relative celerity $ds/dt = \omega$. The canonic form of that wave in non-redundant coordinates is thus $W_t = \cos[(\omega-v) t]$ or, converserly, $W_s = \cos[ (1-v/\omega) s]$.

By taking the definition of $\star$ given in the Appendix \ref{AB} [equation \eqref{star}], the retrieval operation in equation \eqref{op} can be written in terms of contraction\footnote
{The contraction $\odot$ is an \textit{interior product} when a tensor operates onto a differential form (such as $\mathbf e_l$, which is essentially a gradient). When the aforementioned differential form has rank one, the interior product corresponds to a traditional inner product.}
 ($\odot$) and outer ($\otimes$) products as:
\begin{equation}
\mathbf x_{l} = \mathbf X_{s,t} \odot \left( \mathbf e_l^\bot \otimes \mathbf C_{s,t} \right)
\end{equation}
where $\mathbf e_l^\bot$ is the orthogonal complement of $\mathbf e_l$, depending on $l$ alone. 

When space and time are independent, $\mathbf e_s^\bot = \mathbf e_t$ and vice-versa. Here, however, space and time are related by the Spatiotemporal Coevolution Manifold $\mathbf C_{s,t}$, the presence of which means that space and time are not orthogonal complements of each other in the functional sense.

Essentially, the operator performs a contraction or dimensionality reduction in the functional manifold $\mathbf X$ into a submanifold in a defined subspace, say space or time. In practice, $\mathbf x_s$ and $\mathbf x_t$ correspond, respectively, to the spatial and temporal structures of the spatiotemporal process $\mathbf X_{s,t}$.

Rebuilding the spatiotemporal process $\mathbf X_{s,t}$ from its spatial and temporal submanifolds can then be done by introducing an \textit{Interacting Subspace Composition Operator} $\Omega(.)$ and an associated \textit{Composition Product} $\#$:
\begin{equation}\label{galo}
\mathbf X_{s,t} = \Omega(\mathbf x_s , \mathbf x_t) \equiv
\mathbf x_s \, \# \, \mathbf x_t
\end{equation}

The general mathematical definitions of the composition operator and product are introduced in the Appendix \ref{AC}. Their direct application to the spatiotemporal composition in equation \eqref{galo} is given by:
\begin{equation}
\mathbf x_s \, \# \, \mathbf x_t \equiv \mathbf x_s \otimes \mathbf x_t \odot \mathbf C_{s,t}
\end{equation}
where $\mathbf C_{s,t}$, the aforementioned coevolution manifold, defines the rank-$c$ subspace of spatiotemporal interactions, i.e. where $\mathbf x_s$ and $\mathbf x_t$ interact. In practice, it also refers to the relative rate of structural regional and climatological changes, i.e. in space and time, reflecting the dynamic spatiotemporal codependence in the process $\mathbf X_{s,t}$.

All in all, $\mathbf X_{s,t}$ is retrieved from its spatial and temporal subspaces $\mathbf x_s$ and $\mathbf x_t$ by taking their rank-additive outer product and contracting with their interacting subspace spun by $\mathbf C_{s,t}$. The total rank or dimensionality of $\mathbf X_{s,t}$ is thus that of space plus time minus the rank of the space-time interactions.

Unlike subspace decompositions based solely on projections, in our formulation the dynamics are stored in the coevolution manifold, enabling the recomposition of the spatiotemporal process from its spatial and temporal structures. In a non-coevolutionary setting, wherein dim$(\mathbf C_{s,t})=0$, the decomposition and composition operators reduce to the usual projection and composition operations involving linearly separable subspaces. \\

\section{Evaluation of Dynamical Predictability}\label{s4}

\subsection{Dynamic Interaction Analysis}\label{s41}

Having retrieved with DSA a functional basis for the dynamical system of interest, it can now be used to evaluate its predictive power with regards to other processes in which dynamics the system may play an active part.

In this regard, simple linear or nonlinear correlation measures do not suffice, as they provide an aggregate statistical relation without any considerations of dynamic nature. Moreover, they leave out the information content stemming from higher-order nonlinear dynamical interactions.

In order to capture the dynamics at play between processes, we hereby introduce a Dynamic Interaction Analysis (DIA). The procedure basically consists of evaluating the linear and nonlinear dynamical sensitivities between processes, i.e. their dynamic interactions, along with their evolution in time.
This will be a fundamental building block for dynamic model design (section \ref{s42}).

Formally, \textit{the Dynamic Interaction of $k^{th}$ order among $n$ spatiotemporal multivariate processes} $\mathbf M_1,\cdots, \mathbf M_n$, is hereby defined as:
\begin{equation}\label{dip}
\mathcal D^k(\mathbf M_1,\cdots,\mathbf M_n) \equiv \nabla^k_{\left({\mathbf M}_1,\cdots, {\mathbf M}_n\right)} \left(\dot{\mathbf M}_1,\cdots, \dot{\mathbf M}_n\right)
\end{equation}
This generalises Equation \eqref{eq2} for multiple multivariate processes. 
In particular, the first order dynamic interaction between two processes living in a vector space is given by the Jacobian matrix of their coupled dynamical system.

By evaluating all self and cross sensitivities of various orders within the system, the measure in equation  \eqref{dip} captures linear and nonlinear dynamic interactions required to characterise the dynamic feedback structure of the system.

As with the Dynamic Sources, we are also interested in the spatial and temporal structures of the dynamic interactions between processes.

The spatial structure of the interaction of $k^{th}$ order between processes $\mathbf M_i$, $i \in \{1,\cdots,n\}$ captured over a time period $\mathcal T$ can then be expressed by the spatial submanifold of their dynamic interactions, given for each order $k \in \mathbb N$ by:
\begin{equation}\label{dcd_space}
 {\left\llbracket \mathcal D^k \right\rrbracket}_{\mathcal T}
\equiv
\Upsilon_s[\mathcal D^k(\mathbf M_1,\cdots,\mathbf M_n)] 
\end{equation}
where $\Upsilon(.)$ is as in equation \eqref{op} and its operand given by equation \eqref{dip}.

Analogously for time, the temporal structure captured over a certain spatial region $\mathcal S$ is given by:
\begin{equation}\label{dcd_time}
{\left\llbracket \mathcal D^k \right\rrbracket}_{\mathcal S}
\equiv
\Upsilon_t[\mathcal D^k(\mathbf M_1,\cdots,\mathbf M_n)]
\end{equation}

The Spatiotemporal Dynamic Interaction of $k^{th}$ order among the processes $\mathbf M_k$ is then retrievable from their spatial and temporal structures with the $\#$-product introduced in equation \eqref{galo}:
\begin{equation}
\mathcal D^k =
{\left\llbracket \mathcal D^k \right\rrbracket}_{\mathcal T} \, \# \,
{\left\llbracket \mathcal D^k \right\rrbracket}_{\mathcal S}
\end{equation}

\subsection{Dynamic Model Design}\label{s42}

Having retrieved independent dynamic sources $\mathbf X$ (e.g. atmospheric controls) underneath a coevolutionary dynamical system $\mathbf{\dot Y}$ (e.g. geopotential height fields), the goal of the present section is to establish a framework for the simulation of predictand processes $\mathbf Z$ (e.g. precipitation) controlled by $\mathbf{\dot Y}$. This can be addressed as a simple nonlinear dynamic model, whereby the processes of interest $\mathbf Z$ are written as a dynamic function of the sources $\mathbf X$:

\begin{equation}\label{e5}
{\bf \dot Z} = \sum_{k=0}^q k!^{-1} \left[ {\mathcal D^k(\mathbf Z,\mathbf X)} \right]_\mathcal{R} \, {\bf X}^k +  \mathcal O(\mathbf X ^{q+1})
\end{equation}
where $\mathcal R$ denotes evaluation at the reference manifold, the dynamic sources $\bf X$ are governed by equation \eqref{e4}, and $\mathcal O(\mathbf X ^{q+1})$ refers to higher-order dynamics beyond the explicit analytical truncation at order $q$. The $k^\text{th}$ order powers ${\bf X}^k$ account for nonlinear interactions of $k^\text{th}$ order within each component of $\bf X$ but not between different components or sources due to the disentanglement condition in equation \eqref{eq2}.

Equation \eqref{e5} provides a general model structure template or recipe to address the dynamical evolution of any system in terms of a functional basis comprised of underlying independent controls.
\textbf{The determination of a functional basis in a dynamical system is an important step to its optimal modelling, as it ensures that the system is expressed in terms of the smallest number of dimensions without loss of generality.} 
By taking into account not only linear elasticities ($\mathcal D^k$ for $k=1$) but also nonlinear codependencies ($\mathcal D^k$ for $k>1$), this equation can then address changes in dynamical tendencies in terms of nonlinearly interacting processes.

\section{Dynamic Source Analysis of Nonlinear Geophysical Controls}\label{s5}

\subsection{Characterisation of the Geophysical Observable and Associated Datasets}\label{s51}

Our quest for nonlinearly interacting controls in the hydroclimatic system begins with one of the most popular geophysical observables that embodies the dynamical imprint of such interacting processes: the spatiotemporal Geopotential Height fields at pressure levels $p$, $\Phi_p$.

In a nutshell, $\Phi_p$ provides information on the geospatial dynamics of atmospheric waves at the pressure level $p$ to which $\Phi_p$ refers, e.g. $p=500$ hPa. High and low pressure centres can be detected in geopotential terms as high and low heights of a given pressure level. In particular, the large-scale fields provide information on low-frequency atmospheric waves at the planetary scale.

In their natural form, $\Phi_p$ contains not only information on the state of the atmosphere at the pressure level $p$ but also on how it is being affected by underlying forcing mechanisms taking part in its dynamics.
In fact, embedded at the fields for each pressure level lie dynamical footprints of processes at other levels, through vertical interactions that ultimately bring in features stemming from other geophysical domains. 
These include oceanic and land-surface processes at large scales, resulting in geopotential heights actually responding to the difference between the state of the atmosphere over large water and land bodies such as oceans and continents, even considering mid-upper pressure levels such as 500 hPa. Higher pressure levels (corresponding to lower vertical levels) will naturally have a finer sensitivity to the surface properties, at the expense of getting lost in the details and losing the grasp on the bigger picture.

The geophysical mechanisms captured in the geopotential height fields are in turn at the source of the atmospheric circulation regimes that ultimately shape the weather, including moisture advection and thermodynamics, underlying the precipitation dynamics we intend to understand. Comprehensive, introductory treatments on associated atmospheric physical processes are found at \cite{Salby2012} and \cite{Holton2004}.

One of the most reliable and popular datasets with geopotential height fields stems from the NCEP-NCAR Reanalysis Project \citep{Kalnayetal1996}. Therefore, it is the geopotential dataset of choice for our dynamical source analysis.
In the present study, the dynamic source analysis methodology presented in section \ref{s2} is implemented onto NCEP-NCAR Reanalysis data on the 500 hPa Geopotential Height field for a subset of the Northern Hemisphere spanning [30,70] degrees of latitude North and [-80,40] degrees East of longitude. This area covers most of Europe, the Mediterranean and the North Atlantic.

Our interest is in dominant atmospheric processes at the regional and monthly to multidecadal scales. Therefore, we consider datasets with an adequately coarse spatial resolution ($2.5\times2.5$ degree grid-cells), and 30-day moving average of Geopotential Height daily datasets, over a multidecadal period ranging from 1951 until 2012.

\subsection{Statistical Physics of the Geophysical Controls}\label{s52}

In order to evaluate whether traditional statistical tools such as PCA, optimal only for normally distributed linearly mixed data, would be appropriate to the analysis of Geopotential Height fields, we diagnose the information-theoretical measure of Negentropy \cite{Comon1994, Hyvarinenetal2001}. This is a non-negative measure, given by the information entropy deficit of a generic multivariate signal $\mathbf S$ relative to that of a multivariate normally distributed random variable $\mathbf S_N$ with the same mean and standard deviation as $\mathbf S$, i.e.:
\begin{equation}
J(\mathbf S) = H(\mathbf S_N) - H(\mathbf S) \,,
\end{equation}
where $H$ denotes information entropy and $J$ negentropy in the information sense \citep{Shannon1948}. Higher values of $J$ correspond to lower statistical dependence among the components of $\mathbf S$.

In the present study, Negentropy is computed with the recently developed nonlinear statistical estimator of mutual information by anamorphosis \cite{PiresPerdigao2012, PiresPerdigao2013} for better robustness against outliers and data sparsity, noting that information entropy is the mutual information between a variable and itself. The procedure consists essentially on performing a bijective, information-invariant, differential-geometric transformation between the original state space and a working space for optimal statistical analysis without loss of information. Recent applications of the aforementioned method include the estimation of nonlinear information codependences among coevolving landscape-climate processes \citep{PerdigaoBloeschl2014} and the development of an information-theoretic framework for nonlinear scale interactions among chaotic processes in complex systems \citep{PiresPerdigao2015}.

\begin{figure}
\centering
\noindent\includegraphics[width=40.5pc]{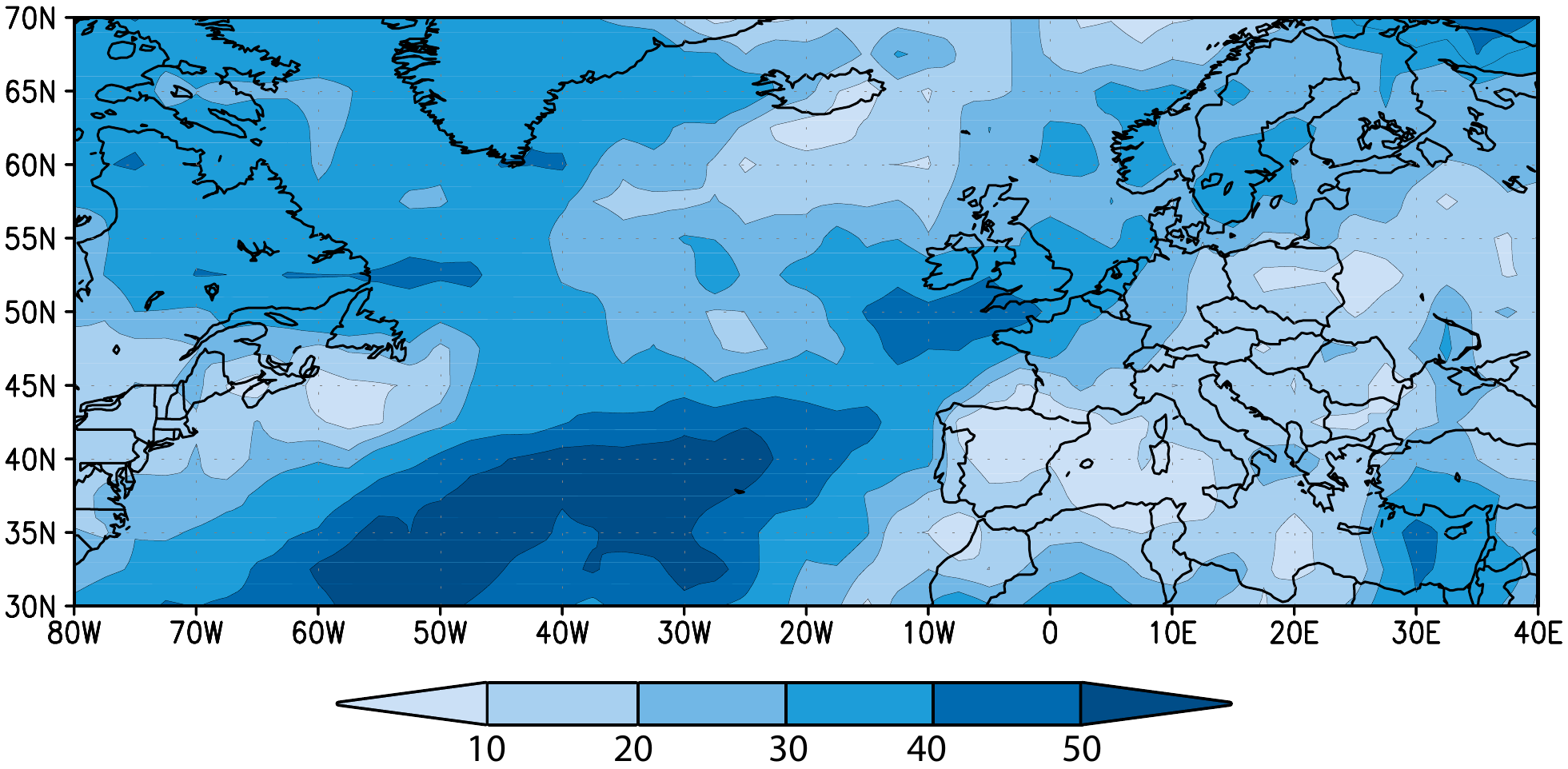}
\vspace{-10pt}
\caption{Negentropy ($mnat$) of the Geopotential Height Fields at 500 hPa considering the joint spatial structure of all climatology partitions within the 1951-2012 period.}
\vspace{10pt}
\label{F3}
\end{figure}

Figure \ref{F3} shows the Negentropy of the 500 hPa Geopotential Height field, with a significance level at 95\% of 0.01 nat (10 mnat). Nat is the unit in which information entropy is expressed when natural logarithms are taken in its computation. For most of the covered area Negentropy is higher than 10 mnat, which indicates a significant departure from Normality, especially over the Atlantic, where Negentropy reaches its maxima. 

Clearly, the Geopotential Height fields are not normally distributed, which indicates that statistical tools assuming or designed for such distribution are not an optimal way to address these fields. 

Naturally though, it might be argued that data could be pre-whitened. However, that would remove precisely the kind of information that is most interesting in the Geopotential Height fields, namely that it is actually not a "free" variable but rather one under the influence of external processes. In fact, by exhibiting negentropy, their distribution exhibits statistical signatures of underlying controls. This statement stems directly from the definition of Negentropy as entropy reduction with respect to the unbounded entropy maximum corresponding to the Normal distribution with the same mean and variance. 
The definition of negentropy further dissipates any misconception about information Entropy as being akin to variance: it is not, as distributions with exactly the same variance can still have completely different information entropies, as evaluated in the negentropy measure.

Higher negentropies diagnose higher entropy reductions, which suggest the existence of stronger controls being enforced on the geopotential over the areas where that reduction is strongest. The geopotential thus serves as proxy for such underlying controls. For instance, while it affects wind patterns, it is also affected by them, which can be understood by bearing in mind that just as a gradient in potential energy drives flow, it is progressively depleted by it as the system works towards geopotential balance. Moreover, albeit the evaluation at a certain pressure level, the horizontal heterogeneity of the geopotential height negentropy informs on entropy constraints imposed by vertical dynamics, including vertical moisture and heat exchanges.

Note that negentropy is higher over the subtropical Atlantic than over Europe, suggesting a stronger oceanic contribution towards constraining the Geopotential even at levels where the earth surface would no longer be expected to play a significant role. The aforementioned contrast is more visible particularly at lower latitudes, where the sea is exerting more influence via enthalpy forcing of the atmospheric pressure gradients and subsequent circulation. This may actually bestow more predictability onto the geopotential, as it is endowed with signatures from oceanic processes with longer memory and driving potential than short-term footprints from intrinsic atmospheric dynamics. Further investigation on these matters is deferred to a subsequent study.

The evidence of regions with a significant non-Normal signature and the associated statistical physics elicit the importance of treating the Geopotential height fields with methodologies that take into account the non-normality of the associated atmospheric processes. Failure to do so may result in undesired results, namely the extraction of features that are not fully independent, and a suboptimal dimensionality reduction of the problem.

\subsection{Dynamic Source Analysis of the Geopotential Height Field}\label{s53}

The theory introduced in sections \ref{s2} and \ref{s3} is now applied to the Geopotential Height Field datasets characterised in section \ref{s51}. A set of spatiotemporal dynamic sources $\bf X$ is extracted from the coevolutionary observables $\bf Y$ (the spatiotemporal geophysical datasets characterised in section \ref{s51}), following section \ref{s24} and under the physical consistency and nonlinear disambiguation constraints in section \ref{s25}.

In practice, the dynamical interactions are evaluated for $k \leqslant \beta$, where $\beta$ is the order of the Sobolev space $\mathcal H^\beta$ defined in section \ref{s2}. This corresponds to the maximum order of weak differentiability of the observable fields, which in the illustrative case considered in the present section yields $\beta=5$.

In the implementation of the DSA procedure, time derivatives up to fifth order are estimated from the spatiotemporal fields with \cite{Geiser2007}'s high-order time discretisation scheme based on the application of the Richardson extrapolation \citep{Descombes2001} to the second-order Crank-Nicolson method \citep{StoerandBurlish2002}. In this regard, the time step is the temporal resolution of the spatiotemporal dataset and the field values for each involved cell in the scheme correspond to those from grid point data.
An analogous reasoning is devised for spatial derivatives quantifying cross-variable sensitivities: in this case, the discrete spatial derivatives are computed using high-order compact ADI (alternating direction implicit) schemes \citep{Karaa2006}. 

While these schemes are originally intended to evaluate discrete derivatives approximating continuous ones (e.g. in computational fluid dynamics), they are also useful for evaluating discrete differentials from datasets, in which case the grid point neighbourhoods provide the scheme inputs natively. The main reason for the aforementioned approach rather than basic discrete derivative schemes stem from the inadequacy of the latter to reliably estimate high-order derivatives from data. 

Without a priori restriction in the number of sources, an objective cut-off criterium is established whereby only dynamic sources sharing effective interaction information with the observables are retained. This leads us to two effective dynamic sources, henceforth denoted \textit{Meridional Dynamic Source} (MDS) and \textit{Zonal Dynamic Source} (ZDS) due to associated properties unveiled below. 

\begin{figure}
\noindent \includegraphics[width=41pc]{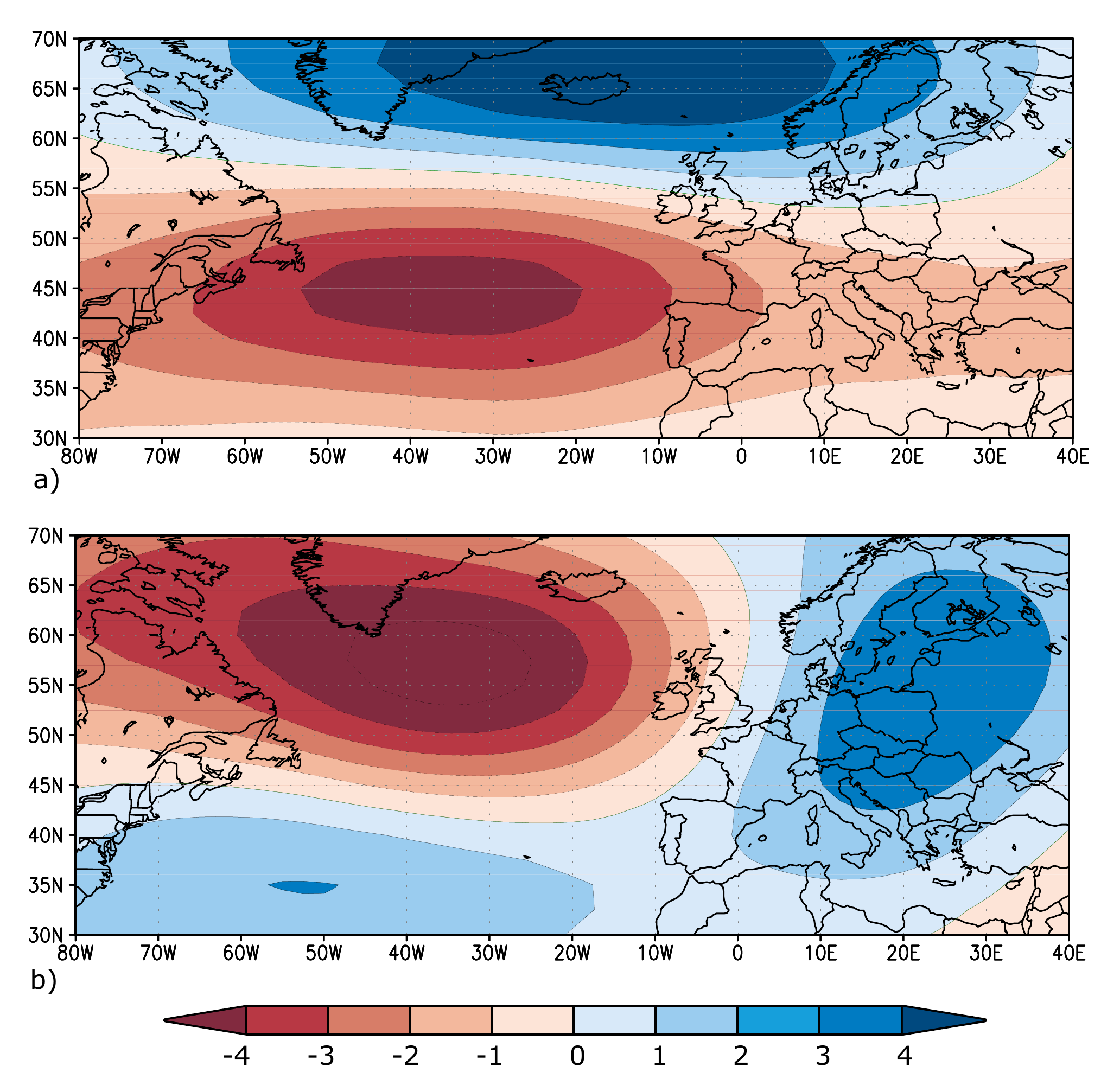}
\vspace{-10pt}
\caption{Spatial Structure of the Meridional (a) and Zonal (b) Dynamic Sources (Functional Basis) of the Geopotential Height Fields at 500 hPa ($X_{\phi}$ and $X_{\lambda}$ respectively), considering the joint spatial structure of all climatology partitions within the 1951-2012 period.}
\label{F4}
\end{figure}

Figure \ref{F4} depicts the spatial structure of the Meridional (Figure \ref{F4}a) and the Zonal (Figure \ref{F4}b) Dynamic Sources $\bf X$ retrieved from the Geopotential Height datasets, i.e. from the observable $\bf Y$ through Dynamic Source Analysis (section \ref{s2}). The spatial structure is computed from \eqref{op} with $l \equiv s$ on $\mathbf X \equiv \mathbf X_{s,t}$, by taking into account all climatology partitions within the period 1951-2012, i.e. of all sub periods obtained by differential partition of the original period. 

The spatial structure of the first Dynamic Source shows an improved geospatial signature of the North Atlantic Oscillation, as compared to spatial patterns obtained by PCA (e.g. \cite{BarnstonandLivezey1987}). In fact, PCA-retrieved geospatial patterns in the literature are uncorrelated but not truly independent, due to the non-normality of the Geopotential Height fields as explained in section \ref{s52}. Such PCA-retrieved patterns are actually different combinations of a smaller set of fundamentally independent processes, which leads to an excessive number of Principal Components and Empirical Orthogonal Functions, when a smaller set of fundamental features suffices to fully characterise the datasets. 

As for the second Dynamic Source of the Geopotential Heights at 500 hPa, its spatial structure suggests the presence of an oscillation pattern with an East-West (zonal) dipole and centres of action over Greenland and the Baltic. Therefore, this is essentially a Zonal Dynamic Source representing a zonal dipole, whose statistical physics are consistent with the Baltic-Greenland Oscillation (BGO) unveiled by \cite{Perdigao2004}. 

Overall, the depiction of both retrieved spatial structures is consistent with the geospatial patterns extracted by \cite{Perdigao2004} with ICA, which had been performed by information-theoretical minimisation of mutual information and maximisation of information negentropy. Relative to ICA, DSA has the added value that not only the statistical climatological aggregate is independent, but also the dynamics within. That way, dynamic information is retained that can be used not only in downscaling but also in dynamic prediction of evolving precipitation distributions (section \ref{s6}).

A common misconception associates statistical patterns with physical regimes. By statistically lumping information, the dynamical sequences upon which regimes should be built are lost. Different dynamical processes can lead to the exact same statistical distribution. That way, while the dynamics explain observed statistics, statistics do not elicit underlying dynamics.

In DSA the spatial structure of the underlying dynamics expresses state-spatial energy optima of the Geopotential Height Fields, i.e. its non-permanent, meta-stable states. 
It should be noted that no perennial stable states exist in the system, as the Earth is permanently far from equilibrium.

\begin{figure}
\centering
\noindent\includegraphics[width=41pc]{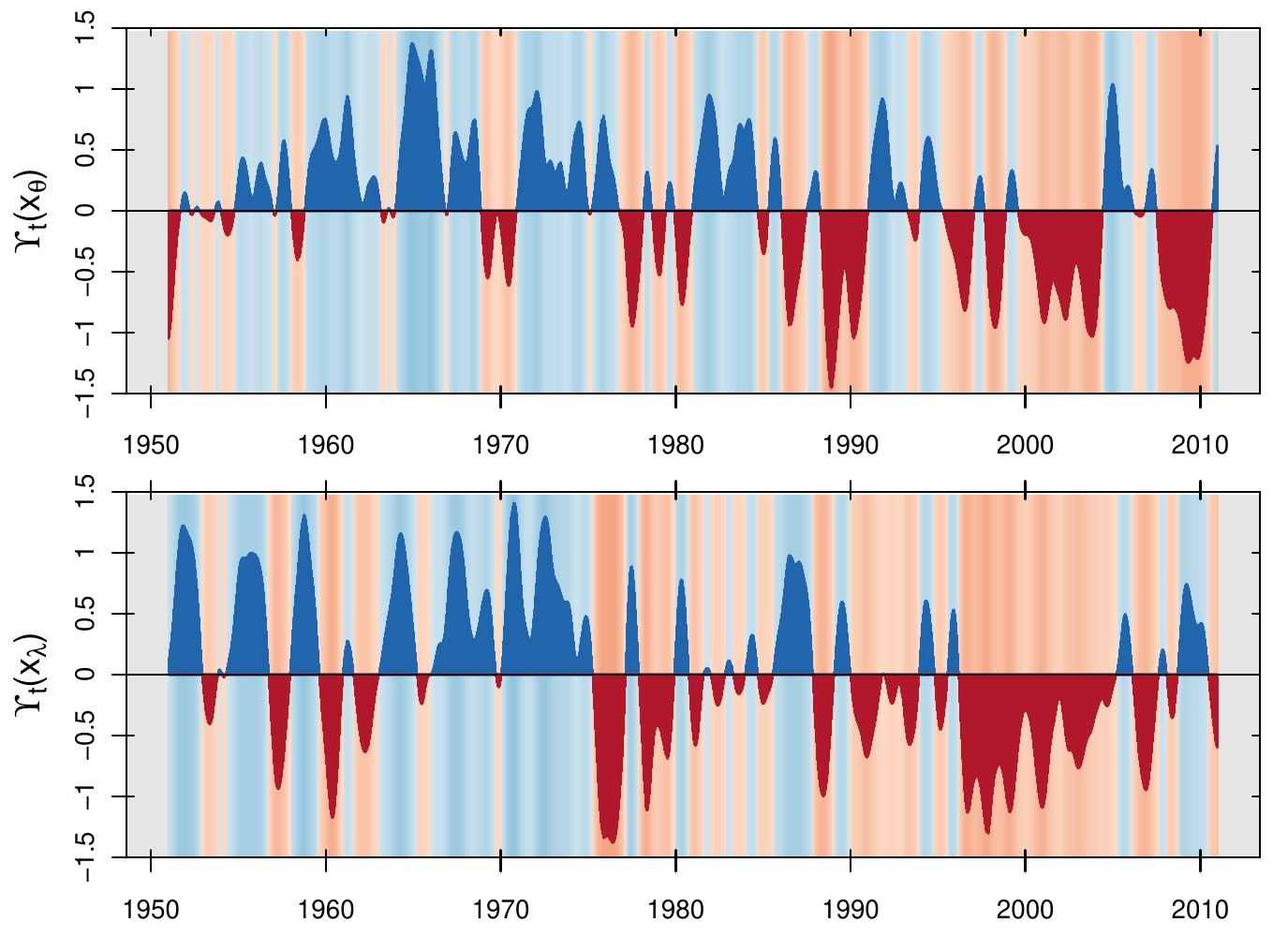}
\vspace{-10pt}
\caption{Temporal (Annual) structure of the (a) Meridional and (b) Zonal Dynamic Sources $X_{\theta}$ and $X_{\lambda}$, considering the joint temporal structure of all regional partitions within the study area.}
\label{F5}
\end{figure}

The temporal structure of the Dynamic Sources is depicted in Figure \ref{F5}, by taking the joint temporal structure of all regional partitions within the study area. The Meridional Source is depicted in (a) and the Zonal in (b). These have been obtained by the temporal subspace retrieval of the respective Dynamic Sources, using the retrieval operation introduced in equation \eqref{dcd_time}.
For visualisation purposes alone, the depicted temporal structures are annually aggregated, which introduces artificial correlations. For that reason, their depiction may appear correlated when in reality the raw dynamic sources on a monthly basis are not.

The values of the temporal structures represent the standardised departure of the respective dynamic sources from their climatological centre of mass. In practice, positive [negative] values in the meridional and zonal cases correspond respectively to northward [southward] and eastward [westward] shifts in the atmospheric wave train and on its associated synoptic systems.

Together with the spatial structure, the temporal structure provides the dynamical evolution of the full geopotential height fields, i.e. their spatiotemporal structure. In mathematical terms, the spatiotemporal structure is the result of the interacting subspace composition (Appendix \ref{AC}) between the temporal and spatial structures.

In the present illustrative case, the manifold relating spatial and temporal structures of the dynamic sources underlying the geopotential height field at 500 hPa is given by the projection of the full spatiotemporal sources onto the coevolution manifold: $\langle \mathbf X_{s,t} , \mathbf C_{s,t} \rangle = \left[ \Upsilon_s(\mathbf X_{s,t}) \odot \mathbf C_{s,t} \right] \otimes \left[ \Upsilon_t(\mathbf X_{s,t}) \odot \mathbf C_{s,t} \right]$. 
From within the multi-order spatiotemporal dependencies, in this particular case the spatiotemporal codependence is only manifested at orders $k>2$. By integrating $\langle \mathbf X_{s,t} , \mathbf C_{s,t} \rangle$ in time over the climatological period and evaluating the spatiotemporal interactions of up to fourth order, a spatial pattern emerges practically coinciding with the Negentropy map in Figure \ref{F3}. Its visual depiction is thus redundant. 

By noting that negentropy captures higher-order statistics, it might not be unreasonable to conjecture that beneath those higher-order statistics lie nonlinear spatiotemporal interactions. These in turn might be related to slow-fast climate dynamic interactions bearing in mind the role of spatial structures in representing processes deemed "slow" relative to the time scale at which the temporal structures are evaluated \citep{PerdigaoBloeschl2014}. 
In particular, these signatures might suggest a nonlinear interplay between "slow" multidecadal and "fast" monthly dynamics in the climate system, with the former conditioning the regimes of the latter, and the latter feeding back on the structural properties of the former. 
A detailed investigation of the dynamical links beneath Negentropy beyond these conjectures makes for a fertile discussion and analysis that is outside of the scope of this paper, rather being deferred to a subsequent study. 

From a fundamental physical standpoint, we interpret meridional dynamic sources in atmospheric dynamics as being associated to the latitudinal (meridional) curvature of the earth with associated differential surface heating by the sun. Meridional thermal gradients between the equator and the poles then result in poleward heat redistribution, which ultimately results in the meridional components of atmospheric flow. If the earth had cylindrical symmetry around its north-south axis, no such meridional dynamics would take place. 

As far as the zonal dynamic sources are concerned, our interpretation consists on the fundamental mechanism underneath being the earth's rotation. This is what fundamentally introduces the longitudinal (zonal) component in the atmospheric dynamics unleashed by the aforementioned meridional gradient. Zonal thermal gradients and associated baroclinicity between large ocean and land masses, while important, come as modulators of the large-scale zonal dynamics induced by planetary rotation. 

Therefore, while observable zonal and meridional flow mutually interfere in the geophysical continuum, the fundamental zonal and meridional processes underneath the dynamics are in essence independent: the meridional gradients come from differential heating due to the shape of the earth, and the zonal from planetary rotation and longitudinal heterogeneities in the energy budget.

\section{Dynamic Predictability and Simulation of Precipitation from the Geophysical Controls}\label{s6}

The quest for geophysical controls on hydro-climatic processes aims at providing a better understanding of fundamental mechanisms and their relevance to the ultimate behaviour of the hydroclimatic system. Having retrieved and discussed such controls, the aim of the present section is to relate them with precipitation regimes: firstly by diagnosing dynamic relations between them, and secondly by dynamically predicting precipitation dynamics from the geophysical controls. 

\subsection{Diagnosing Dynamic Predictability}\label{s61}

We begin by analysing measures of dynamic dependence between Precipitation and the Dynamic Sources of the Geophysical Controls for all climatology partitions within the study period.

Formally, the spatial structure of the dynamic interaction of $k^\text{th}$ order of the dynamic process $\bf Z$ (e.g. Precipitation) relative to the sources $\bf X$ (e.g. MDS, ZDS) over a climatology can be expressed by taking equation (18) with  $\mathbf M_i$ given by $(\mathbf M_1,\mathbf M_2) = (\mathbf X,\mathbf Z)$:
\begin{equation}\label{dcd}
\left \llbracket \mathcal D^k(\mathbf X,\mathbf Z) \right \rrbracket_{\mathcal T} \equiv \Upsilon_s \left[ \mathcal{D}^k_{\mathcal T}(\mathbf X,\mathbf Z) \right]
\end{equation}
where $\Upsilon_s(.)$ is as in equation \eqref{op} with subspace $l\equiv s$, $\llbracket \cdot \rrbracket_\mathcal T$ denotes, in this application, the joint spatial structure of all climatology partitions within the 1951-2012 study period.

In practice, a normalised dynamic codependence is computed from equation \eqref{dcd} as follows: 
\begin{equation}\label{kpred}
\mathcal N^k_{\mathcal T} (\mathbf X,\mathbf Z) \equiv \left \llbracket \mathcal D^k_{\mathcal T}(\mathbf X,\mathbf Z) \right \rrbracket_{\mathcal T} \, \left \llbracket {\Lambda_k} \right \rrbracket_{\mathcal T}^{-1}
\end{equation}

Equation \eqref{kpred} yields essentially an adimensional measure of $k^{th}$-order predictability among the intervening variables at each point on a map: in particular, \textit{Linear Predictability} for $k=1$ and \textit{added value of Nonlinear Predictability} for $k>1$. The absolute value of the measure ranges from zero in dynamically independent cases to one in fully redundant cases at the respective order. 

An integration of this measure over a sub manifold of the dynamical system yields an overall statistical predictability consistent with the Information-Theoretical Correlation introduced by \cite{PiresPerdigao2007} as a generalisation of the Equivalent Correlation from \cite{Perdigao2004}. In particular for the linear case $(k=1)$, the integration of the measure yields an aggregate diagnostic corresponding to a traditional correlation ranging from -1 to 1.

In practice, we are interested in effective interactions, i.e. the measured dynamic interactions discounting the spurious interactions between a given number NS of independently Monte Carlo (MC) shuffled dynamic sources and precipitation. 

For that purpose, let $\mathcal N^k_{\mathcal T}(\mathbf X_\text{MC},\mathbf Z_\text{MC})_\text{NS}$ denote the ensemble mean of the dynamical interactions between NS randomly shuffled realisations of $\bf X$ and $\bf Z$, with independent shuffling applied to each variable. Then the effective interaction of order $k$ is given by
\begin{equation}
\mathcal N^k_{\mathcal T}(\mathbf X,\mathbf Z)_{\text{eff}} = \mathcal N^k_{\mathcal T}(\mathbf X,\mathbf Z) - \mathcal N^k_{\mathcal T}(\mathbf X_\text{MC},\mathbf Z_\text{MC})_\text{NS}
\end{equation}
In the practical implementation considered in the present section, NS=$10^4$.

The precipitation data has been retrieved from the Global Precipitation Climatology Centre (NCAR, 2014), and contains global analysis datasets of monthly precipitation on the earth's land surface based on \textit{in situ} rain gauge observational data. The data processing procedures conducted by GPCC are detailed in e.g. \cite{Rudolfetal1994}, \cite{RudolfSchneider2005} and involve statistics exclusively over the observational data.

\subsection{Dynamic Predictability Results and Discussion}\label{s62}

\begin{figure}
\noindent\includegraphics[width=40.5pc]{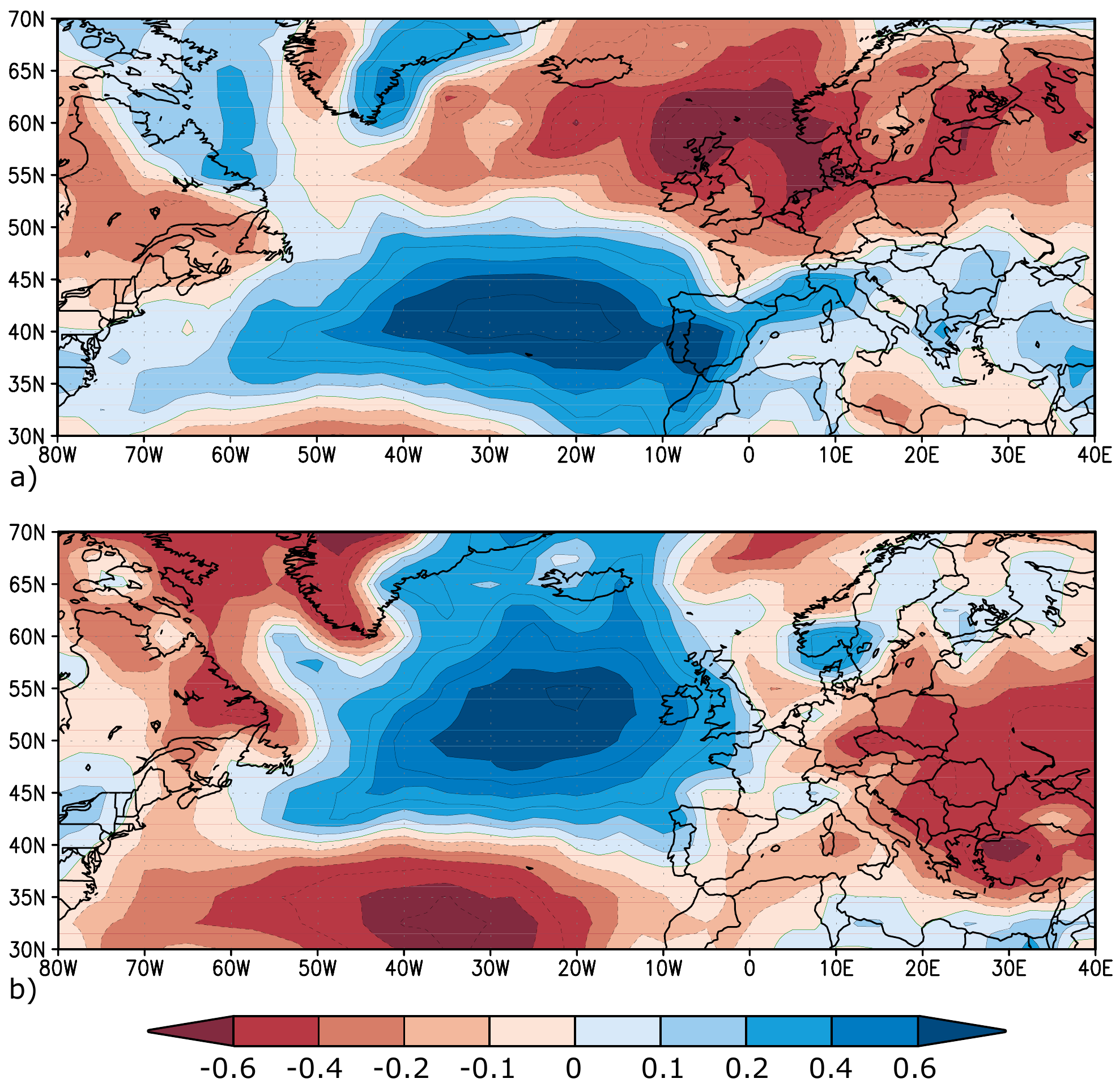}
\vspace{-10pt}
\caption{Linear Predictability of Monthly Precipitation to the Meridional (a) and Zonal (b) Dynamic Sources ($X_{\phi}$ and $X_{\lambda}$) depicted in Figure \ref{F4}, again considering the joint spatial structure of all climatology partitions within the 1951-2012 period.}
\label{F6}
\end{figure}

\begin{figure}
\noindent\includegraphics[width=40.5pc]{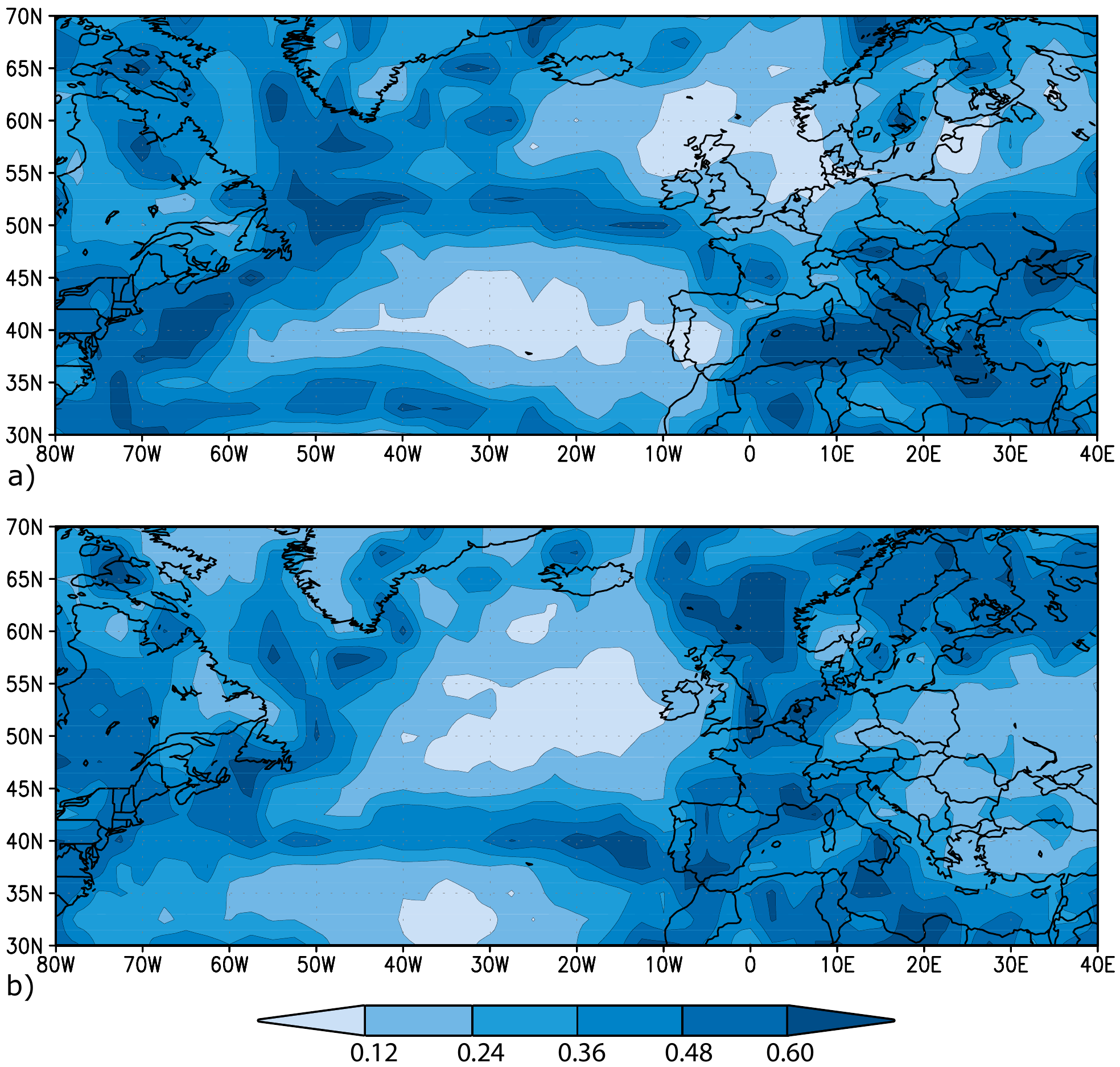}
\vspace{-10pt}
\caption{Added value from Nonlinear Predictability of Monthly Precipitation to (a) $X_{\phi}$ and (b) $X_{\lambda}$, considering the joint spatial structure of all climatology partitions within the 1951-2012 period.}
\label{F7}
\end{figure}

Figures \ref{F6} and \ref{F7} depict, respectively, the spatial structures of the Linear Predictability ($k=1$) and the added value from Nonlinear Predictability ($1<k\leqslant q$ for $q=5$) of the Monthly Precipitation to the Meridional (a) and Zonal (b) Dynamic Sources, depicted in Figure \ref{F4}, considering the joint spatial structure of all climatology partitions within the 1951-2012 period. The computations have been performed from equations (23)-(25) with $\mathbf X$ representing the Dynamic Sources and $\mathbf Z$ the Monthly Precipitation.
A truncation order $q=5$ was established in practice since higher-order terms no longer added significant value to the predictability. 

The linear predictability is quite considerable over the centres of action of the underlying atmospheric oscillations, namely the Meridional Dynamic Source (MDS) in (a) and the Zonal Dynamic Source (ZDS) in (b). However, it degrades with the distance to those centres of action. As far as the MDS is concerned, this is consistent with the findings of statistical studies conducting linear sensitivity analysis of monthly precipitation to related atmospheric teleconnection patterns (e.g. \cite{MurphyandWashington2001} focusing on the Irish and British Isles and \cite{Trigoetal2004} on Iberia), and also with the first-order components of the nonlinear information-theoretical analysis conducted by \cite{PiresPerdigao2007} at the Euro-Atlantic scale. 

Regarding the ZDS, no such studies are known to exist since the oscillation per se is a new finding, providing the dynamical basis underneath the information-theoretical pattern of the BGO unveiled by \cite{Perdigao2004}. In fact, while the literature features teleconnection patterns with a somewhat related zonal component (e.g. Eurasia-2 pattern in \cite{BarnstonandLivezey1987} known as the Eastern Atlantic / Western Russia teleconnection), the Eurasia patterns are neither statistically independent nor dynamically based, rather relying on hierarchies of explained variance in climatological anomalies of geopotential height fields.

The "lost" predictability away from main centres of action seen in Figure \ref{F6} is essentially "recovered" when the nonlinear contributions to predictability are taken into account (Figure \ref{F7}). In fact, the added value is highest at the linear-poorest regions, essentially closing the gap in predictability between the regions with higher linear predictability and those with scarcer one.
As far as the MDS is concerned, this is consistent with the findings from \cite{PiresPerdigao2007} regarding the nonlinear sensitivity of precipitation to the NAO, to which the statistical physics of MDS are closely related. Here, as there, regions where the linear sensitivity is weak (e.g. Central Europe) actually exhibit significant predictability once the information associated to nonlinear relations between precipitation and NAO is taken into account. In the present work this result is thus revisited with the newly extracted Meridional Dynamic Source, now with a nonlinear dynamical perspective on the NAO complementing the previous nonlinear statistical, information-theoretical one, and further extended to the Zonal Dynamic Source capturing the dynamics underlying the statistical physics of BGO from \cite{Perdigao2004}. 

All in all, the influence over precipitation is mostly linear over the centres of action and nonlinear away from them. By taking into account the full nonlinear information, then the overall domain exhibits strong predictability of precipitation from the dynamic sources.

While the linear information refers to the dynamics associated to mean atmospheric flows, the nonlinear stems from higher-order interactions such as turbulent fluxes. These are particularly important along storm tracks and regions dominated by convective behaviour. By bringing out the nonlinear predictability, DSA elicits not only the traditional large-scale synoptic regimes, but also hotspots of convectivity elusive to the linear terms.

The physical reason why smaller convective and turbulent processes can be dynamically inferred from large-scale atmospheric processes resides in the cooperative nature of synergistic interactions formulated in DSA. In fact, primary processes can synergistically cooperate to produce secondary processes at different scales (e.g. large-scale planetary waves resonate into finer-scale secondary perturbations in the atmosphere) via wave resonance across scales, which can be rigorously characterised in both nonlinear dynamic and information-theoretical terms \citep{PiresPerdigao2015}. 

Our dynamical approach complements overall statistics by enabling the extraction of structural information within the domain, namely in its sub-partitions along with its dynamical evolution, i.e. a dynamic rather than lumped approach. 
By not simply averaging spatiotemporal fields in space and time upon separation, the coevolution manifold is preserved. Therefore, the spatiotemporal interplay is retained, enabling space and time to communicate and accounting for structural changes in space, with the added predictability that ensues under non-invariant non-ergodic climate dynamics. 

Moreover, the consideration of all partitions within the spatiotemporal domain of study captures not only the broader Euro-Atlantic (spatial) and 60-year (temporal) scales, but also a whole spectrum within, down to spatiotemporal resolution limits (time step, grid size). Therefore, each pixel in the spatial structures contains not a temporal mean over a single period, but rather an integrated account of dynamic processes at multiple scales.

\subsection{Dynamic Simulation of Monthly Precipitation from the Geophysical Controls}\label{s63}

The dynamical simulation model [equation \eqref{e5}], introduced in section \ref{s42} as a general simulation framework, is now applied to the simulation of Monthly Precipitation from the Dynamic Sources of the Geophysical Controls disentangled in section \ref{s53} with the methodologies devised in section \ref{s2}.
Taking $\mathbf Z = \mathbf P$ and $q=5$ in equation \eqref{e5} bring us to a practical model form:
\begin{equation}\label{ddmodel_prec}
{\bf \dot P} = \sum_{k=0}^5 k!^{-1} \left[ \mathcal N^k_\mathcal{T} ({\bf X},{{\bf P}})_\text{eff} \right]_{\mathcal R} \, {\bf X}^k +  \mathcal O(\mathbf X^{6})
\end{equation}
where $\bf \dot P$ prescribes the simulated dynamics of monthly precipitation, and ${\bf \dot X}$ is given by equation \eqref{e4di} with $\mathbf X = (X_\theta, X_\lambda)$ being the set of dynamic sources (meridional and zonal) of the geopotential height fields at 500 hPa determined in the previous section, and the model is truncated to the $5^{th}$ order in $\mathbf X$.

\textbf{The model addresses the dynamic evolution of precipitation as a function of the dynamic sources representing the dominant modes of spatiotemporal dynamical variability in the geopotential height fields.}
Moreover, given the reciprocal nature of the Dynamic Interactions, it enables retroaction of Precipitation as it responds to the Dynamic Sources. This is physically important since the weather systems work towards depletion of the baroclinicity at their genesis, a negative feedback stemming from the Second Law of Thermodynamics.
While the fundamental basis functionals $\mathbf X$ are independent, their impact on $\mathbf P$ as formulated in $\mathcal N^k_\mathcal T(\cdot)$ is indeed adjusted (internal feedback) as the dynamics unfold in \eqref{ddmodel_prec}.

The dynamic model is initialised in a two-step procedure: 1) evaluating the long-term dynamical properties of the reference manifold ($\mathcal R$), and 2) the short-term spatiotemporal divergence of the phase space flow from $\mathcal R$. The evaluation in 1) is done via the spatial structure of the dynamic interaction between $\mathbf P$ and the atmospheric controls represented by the dynamic sources $\mathbf X$. In this regard, the spatial structure plays the role of the spatial legacy of long-term dynamics as in \cite{PerdigaoBloeschl2014}. The procedure in 2) finds the phase space directions that maximise the equienergetic entropy production, by concentrating the Lyapunov spectrum \citep{Ott2002} along the highest positive values whilst conserving the spectral integral and thus the energy of the system in that initialisation step. This procedure yields similar effects to the dynamical breeding performed in \cite{Perdigao2010} in that it projects the dynamics onto the most unstable directions, thus maximising the entropy production of the system in line with the physical principles. However, it differs from \cite{Perdigao2010} in that the operation is performed directly over the Lyapunov spectrum of the datasets.

Essentially, the initialisation leads to an initial distribution maximising the ability to further produce entropy in a realistic manner. In fact, while the reference manifold in 1) provides the structure of the slow dynamics (the climate or "personality" of the system), the spatiotemporal divergence in 2) provides the dominant lines of deviation from that core, i.e. the seeds for the faster weather deviations from the climate norm ("mood swings"). 
The dynamics then unfold as prescribed by \eqref{ddmodel_prec} and \eqref{e4di}.

\subsection{Dynamic Simulation Results and Discussion}\label{s64}

The dynamic simulation of monthly precipitation from the underlying dynamic sources has been performed for two illustrative regions with predominantly Atlantic and Continental climatic norms: the former roughly encompasses the Loire Region in France, and the latter the Upper Danube Region in Central Europe. In geodesic coordinate terms, the simulations have been performed over the domains $[45.5\degree N,48\degree N]\times[-1.5 \degree E,3 \degree E]$ and $[47.5 \degree N,49 \degree N]\times[14 \degree E,19 \degree E]$, respectively. 
The denomination of these climate regions should not be mistaken for those of underlying hydrologic basins. 

\begin{figure}
\centering
\noindent\includegraphics[width=37pc]{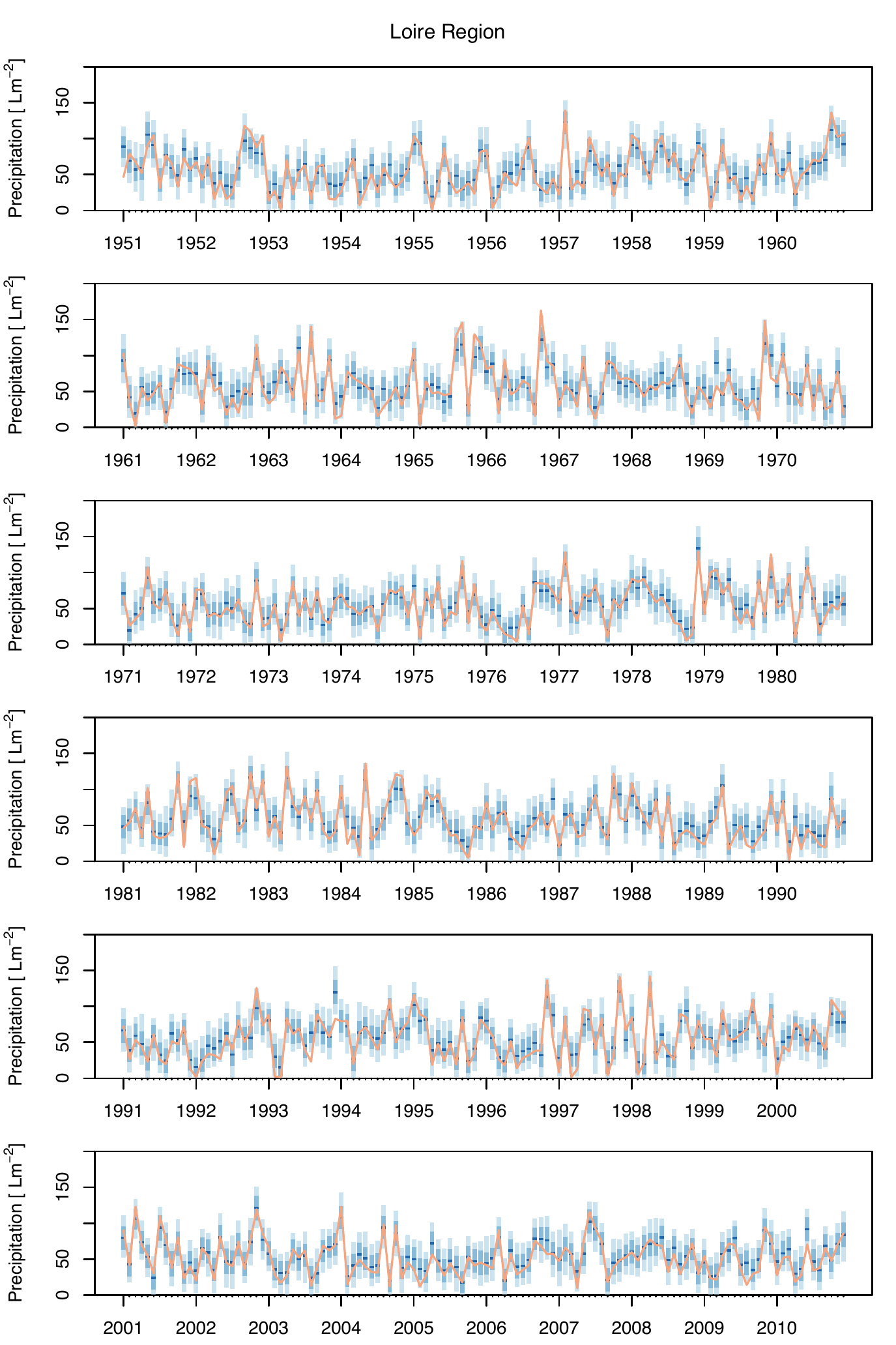}
\caption{Dynamic Simulation of the monthly precipitation from the interacting Zonal and Meridional Sources (blue distributions with median in darker blue) and reference observational time series (orange line) for the Loire region.}
\label{F8}
\end{figure}

\begin{figure}
\centering
\noindent\includegraphics[width=37pc]{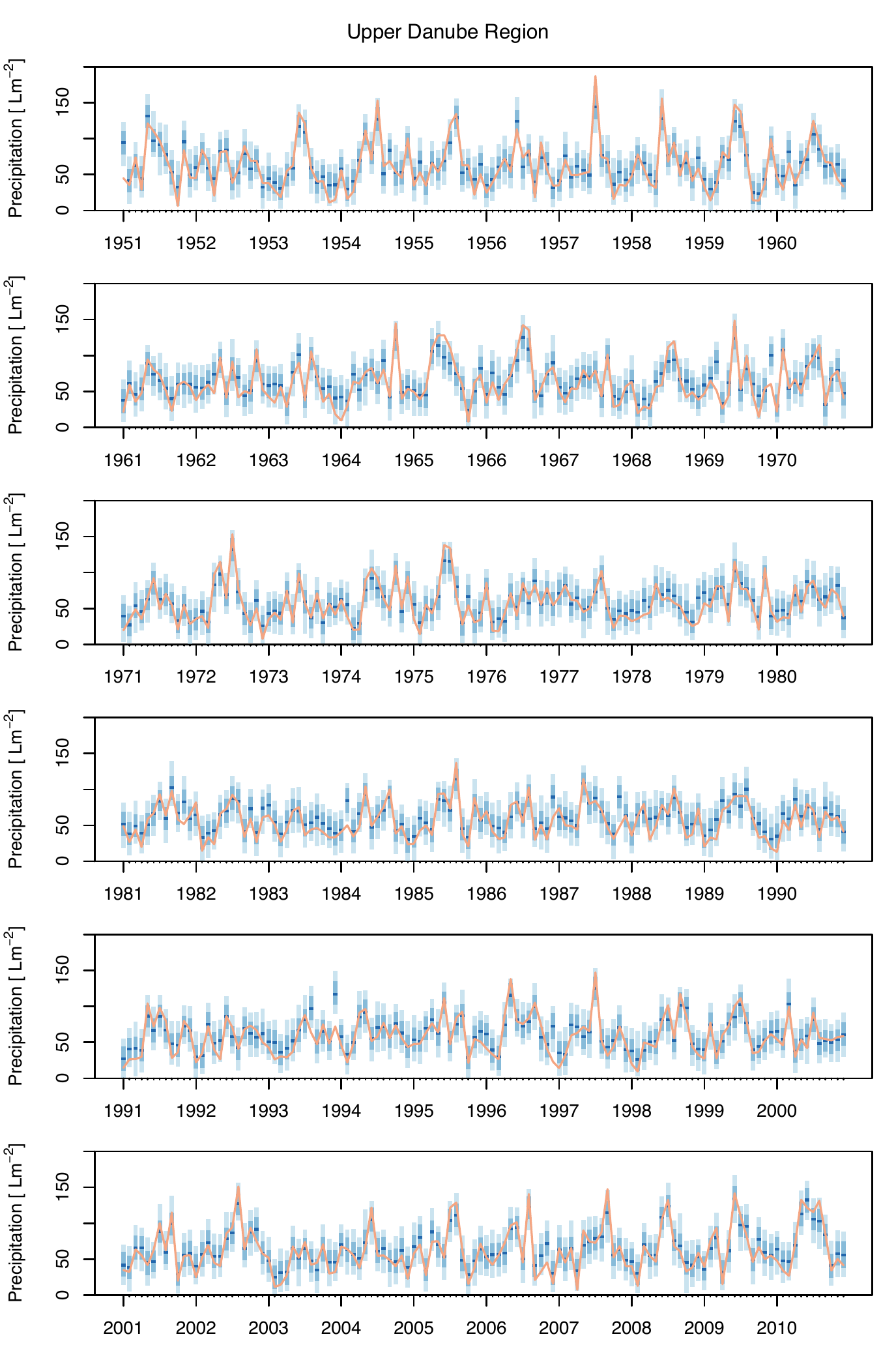}
\caption{Dynamic Simulation of the monthly precipitation from the interacting Zonal and Meridional Sources (blue distributions with median in darker blue) and reference observational time series (orange line) for the Upper Danube region.}
\label{F9}
\end{figure}

\begin{figure}
\noindent\includegraphics[width=41pc]{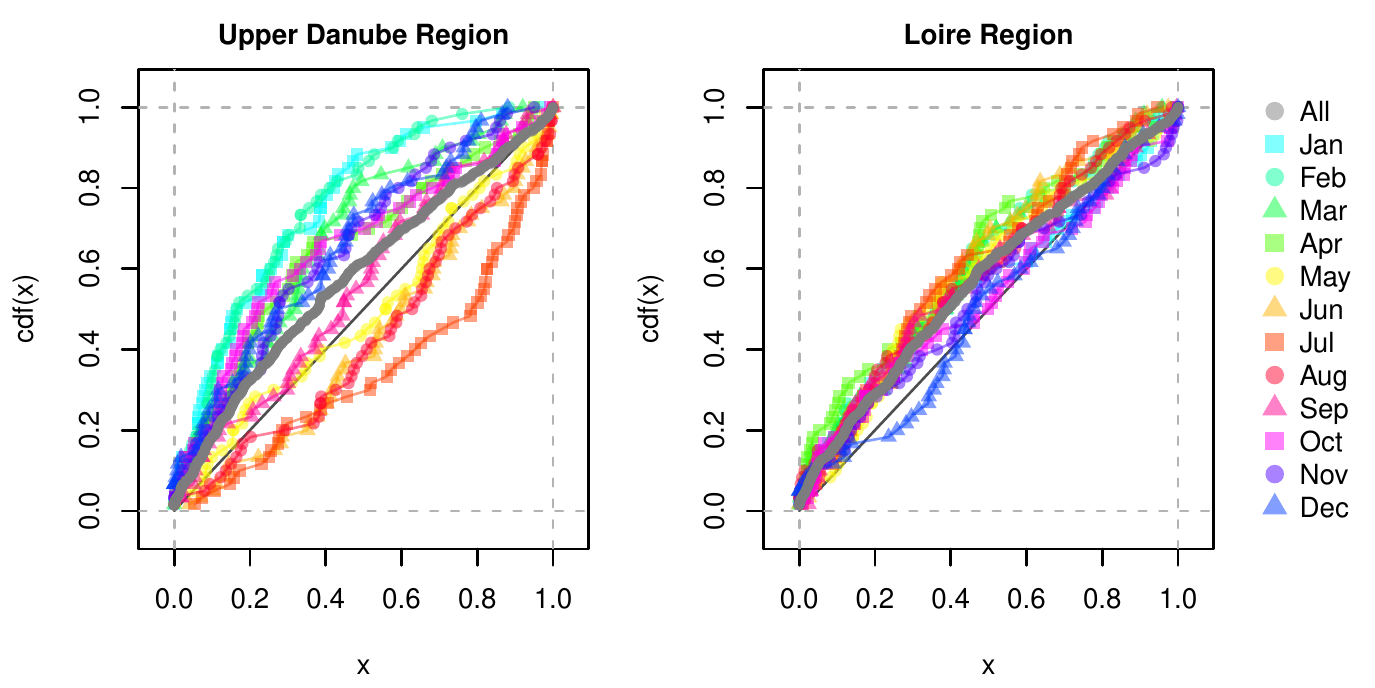}
\vspace{-20pt}
\caption{Quantiles of the simulated precipitation distributions in which the observational series lie. With respect to Figures \ref{F8}-\ref{F9}, the coloured symbols depict the quantile position of the orange lines in the blue-toned distributions.}
\label{F10}
\vspace{10pt}
\end{figure}

Overall, the observed Monthly Precipitation ($P_o$) is captured by the distribution produced by the dynamic simulation model ($P$). This is seen in Figures \ref{F8} and \ref{F9} for the Loire and Upper Danube regions, respectively, where the observations (orange) are found within the modelled distribution (shades of blue), with very few exceptions.
In general, the highest values of $P_o$ sit at the highest quantiles of $P$, consistent with those values being upper extremes in the climatology. 
The converse happens for lower values of $P_o$: in fact, these sit at the lower tail of the $P$ for the corresponding months. 
These considerations are further supported by Figure \ref{F10}, showing in which quantile of $P$ the observation $P_o$ sits. The shape of the cumulative density function (cdf) is due in part to the fifth-order truncation of the model. 

During the drier seasons (Summer in Loire, Winter in Upper Danube), the observed precipitation is captured by lower quantiles of the simulated precipitation, i.e. below their median. Since the outcome of the simulations is not the median but rather the full distribution for each month, and since the distribution does capture the observations at different quantiles, this behaviour does not constitute a modelling bias. 

Whether the observations are captured by upper or lower quantiles of the simulated distributions actually comes down to physical arguments. For instance, higher-order turbulent processes associated with storms are captured by upper quantiles of the simulated precipitation. As for the aforementioned dry season behaviour, the observed precipitation naturally sits in lower quantiles of the simulated distributions due to the relative scarcity of precipitable water for precipitation to occur at the levels that would otherwise be expected given favourable low pressure conditions. 

While moisture content in the atmosphere is not explicitly characterised, its influence is nonlinearly embedded in the geophysical fields. In fact, since moist air has lower density than drier air, it introduces lower pressure anomalies relative to those of drier air conditions. Therefore, air moisture implicitly contributes to the predictability of precipitation from the geophysical dynamics in the atmosphere - a contribution that would be elusive if only the linear information contained in the geopotential height fields would be taken into account.

By simulating an overall distribution, the model captures observations arising from a diversity of situations with associated quantiles as noted above. In doing so, the model captures not only precipitation outcomes associated to first-order atmospheric dynamics (mean atmospheric flows, large-scale synoptic processes) but also regional extremes arising from higher-order processes. 
An ensemble of simulated means (the usual paradigm in ensemble prediction) would not aptly capture the behaviour of extremes, as the physics governing the overall distributions entails laws differing from those governing the first moment. This is why our simulations take full distributions rather than ensembles of simulated means. 

\break

\section{Concluding Remarks}\label{s7}

Dynamic Source Analysis (DSA) was introduced as a novel nonlinear dynamic analysis and modelling framework for improving the fundamental physical understanding and predictability of complex coevolutionary systems. These were formulated in terms of the smallest set of generative fundamental non-redundant processes, i.e. a functional basis to the dynamical systems. In other words, DSA provided a model building and optimisation approach wherein a canonic structure is brought out for expressing the full system complexity in terms of the minimum set of components without loss of information.
This enables data compression beyond the traditional feature extraction and information theoretical methods, by disentangling nonlinearly independent components not only over statistical aggregates but also within each step of the dynamics. 

Given the natural mathematical ambiguity in solutions from nonlinear problems, the uniqueness of the DSA basis is ensured with first principles in Physics, namely on energy and entropy production. For that purpose, the physical solutions are dynamically selected from within the mathematically outcomes by automatically taking the energy minima in the solutions' isentropic subspaces and entropy maxima in their equienergetic subspaces.

The generative ability of the DSA-retrieved functional basis stems from the formulation of synergistic interactions expressing cooperative dynamics among fundamental processes. In doing so, new features are brought as secondary processes that do not exist in any of the intervening primary processes. Further synergistic interactions among primary and secondary processes lead to a third generation of processes and beyond, in a synergistic cascade ultimately leading to the full representation of the system complexity. 

While coevolutionary interactions entail mutual influence and thus redundancy, purely synergistic ones entail only innovative outcomes. A model structure grounded on the latter interactions is thus more efficient in representing complexity in the simplest way without loss of generality. 

In practice, while all coevolutionary redundancies are cleared out between processes, self-interactions are preserved within each process (e.g. entropy production stemming from internal positive feedbacks). This irreducible coevolutionary behaviour is natural, as a process is redundant with itself. The coevolutionary nature of the problem is then concentrated along the "diagonal spine" of the dynamical system, i.e. the diagonal of its dynamic interaction tensor, corresponding to self-interactions within each basis component.  By taking the spectral view on such components, such self-interactions account for internal scale interactions or spectral coevolution, ranging from separable slow-fast dynamics to non-separable interactions among multiple indistinguishable scales (\textit{broadband coevolution}). 

In a spatiotemporal setting wherein the spatial structure embodies the legacy of long-term dynamics, spectral coevolution then comes down to a fundamental dynamic dependence between space and time. In this work, a coevolution manifold was introduced representing this codependence across the overall dynamical system, enabling space-for-time substitution to be performed at once for non-ergodic complex systems exhibiting diverse coevolution rates across different directions in space and time. This opens new perspectives in multivariate space-for-time trading, e.g. in regional hydrology. 

In the hydro-climatic context, DSA was implemented as a means to retrieve a non-redundant set of nonlinearly interacting processes exerting control over precipitation in space and time. For that purpose, DSA has been implemented on the geopotential height fields at 500 hPa, retrieving a nonlinear spatiotemporal functional basis comprising two independent dynamic sources: the Meridional Dynamic Source (MDS) associated to the North Atlantic Oscillation (NAO), and the Zonal Dynamic Source (ZDS) embodying a pressure dipole between the Baltic and Greenland and providing a dynamic basis to the Baltic-Greenland Oscillation (BGO) unveiled by \cite{Perdigao2004}. The latter generalises and dynamically links the blocking systems found over the respective centres of action, namely those over Greenland and over the Baltic-Scandinavian area.

Whilst the NAO is responsible for the North-South shifting of the storm tracks over the Atlantic (with consequences downstream of the atmospheric flow well into Europe), the BGO does so in the West-East shift, counterbalancing or enhancing the zonal flow and thus being responsible for increasing or reducing the blocking risk in synoptic atmospheric processes. In fact, the BGO is associated to the zonal maritime-continental pressure gradients and associated atmospheric circulation, being affected by differential changes in temperature over land vs. over sea (which affects the zonal land-sea baroclinicity and thus the zonal circulation and progression of the synoptic systems).

The newly obtained dynamical sources were decomposed in space and time, by taking into account the spatiotemporal codependence as expressed by relative celerities in the atmospheric dynamics. For that purpose, a new spatiotemporal decomposition methodology of nonlinearly connected manifolds has been introduced and implemented, resulting in the spatial and temporal structures of the dynamical sources.

These have then been used to evaluate linear and nonlinear dynamical interactions and predictability between precipitation and the geopotential heights. The added value of the nonlinear interactions has been put into evidence, showing that by accounting for these the predictability of precipitation from the dynamical sources can be strongly improved throughout the domain.

On the temporal domain, the dynamic simulation of evolving Precipitation distributions has aptly captured the observational values for the test regions, as the latter were mostly captured by the former. Naturally though, higher[lower] observed values mostly lie on higher[lower] quantiles of the simulated distributions for each monthly record, with upper and lower observational extremes being at the corresponding tails of the simulated distributions.

All in all, Dynamic Source Analysis elicits fundamental processes and interactions at play in complex systems, eliminating coevolutionary redundancies and expressing complexity in the simplest dynamic architecture without loss of generality: a synergistic dynamical system. 
Moreover, by providing solutions constrained on first principles in physics, DSA ensures physical consistency in the resulting model structure and intervening processes. 

With these advances in mind, a diversity of applications can thus be considered for the proposed theory of Dynamic Source Analysis, such as signal processing, feature extraction and inverse dynamic modelling from large datasets \textit{(big data)}, complex network design, dynamic model building and architecture optimisation of dynamical systems. 

Further applications to hydro-climate dynamics and analytics are currently in progress, with particular emphasis on improving the dynamical understanding, model design and predictability of coevolutionary regimes, transitions and extremes in hydro-climatic systems.

\appendix

\section{Coevolution vs. Synergies in Information-Theoretical Terms}\label{AA}

The fundamental difference between synergy and coevolution is that while the former involves 
cooperation between independent source components in the production of mutually dependent child processes, the latter entails dependence among the intervening components. In information terms, synergy produces extra information exceeding the sum of the information content of the parts (constructive interference pattern), whereas coevolution entails redundancy, resulting in the total information content being less than the sum of the informations of each individual components.

In order to illustrate these concepts, consider a pair of independent spatiotemporal dynamic sources $\mathbf X_A,\mathbf X_B$ and a set of observable child processes $\mathbf Y$ dynamically forced by the sources, corresponding to a complex dynamical system. The triadic interaction information among them is derived from \cite{PiresPerdigao2015}'s equation (9d), yielding:
\begin{equation}\label{ITPP}
I_T(\mathbf X_A,\mathbf X_B, \mathbf Y) \equiv I[ (\mathbf X_A, \mathbf X_B) , \mathbf Y]  - I(\mathbf X_A, \mathbf Y) - I(\mathbf X_B,\mathbf Y) 
\end{equation}
where $I$ is Mutual Information (always non-negative) and $I_T$ is Interaction Information \citep{Tsujishita1995}.

While $\mathbf X_A$ and $\mathbf X_B$ are independent from each other, they exhibit joint and individual codependence with the child processes $\mathbf Y$. Equation \eqref{ITPP} quantifies the difference between their joint contribution and the sum of their individual contributions to the coevolutionary system represented by $\mathbf Y$. If the former exceeds the latter, $I_T (\mathbf X_A,\mathbf X_B,\mathbf Y) >0$ and we are in the presence of a net synergy. The independence between $\mathbf X_A$ and $\mathbf X_B$ ensures that $I_T$ will at least be non-negative. This can be further seen by expressing the triadic interaction information as the difference between conditional and non-conditional information terms:
\begin{equation}\label{ITPPc}
I_T(\mathbf X_A,\mathbf X_B,\mathbf Y)=I(\mathbf X_A,\mathbf X_B | \mathbf Y) - I(\mathbf X_A,\mathbf X_B)
\end{equation}
In fact, when the sources are independent, $I(\mathbf X_A,\mathbf X_B)=0$ and thus $I_T(\mathbf X_A,\mathbf X_B,\mathbf Y) \geqslant 0$.

Naturally though, if the sources were codependent and their mutual information exceeded the conditional term in the first \textit{rhs} (right-hand-side) term of \eqref{ITPPc}, then the converse could happen, i.e. $I_T<0$ expressing net redundancy.

When looking at the overall statistic-dynamic properties of a complex system, coevolutionary and synergistic interactions can be captured in information-theoretical terms and linked to underlying dynamical systems, as done in \cite{PerdigaoBloeschl2014} and \cite{PiresPerdigao2015}. In our proposed theory of Dynamic Source Analysis, a dynamic treatment is sought wherein coevolution and synergies are expressed in dynamic interaction terms. These aim to capture not only the overall nonlinear information as in the illustrative example \eqref{ITPP}-\eqref{ITPPc}, but also the underlying dynamic relationships that ultimately explain the diagnosed statistics.

\section{Interacting Subspace Retrieval and Decomposition}\label{AB}

\subsection{Retrieval Product}

Consider the generic $n$-dimensional functional or rank-$n$ tensor $\mathbf X$ and the $m$-dimensional functional or rank-$m$ tensor $\mathbf Y$, along with a generic codependence $c$-dimensional functional $\mathbf C$ representing their interdependence.

We hereby define the \textit{Retrieval Product} of $\mathbf Y$-structure from $\mathbf X$ as:

\begin{equation}\label{star}
\mathbf X \star \mathbf Y \equiv \mathbf X \odot \mathbf Y^\bot \otimes \mathbf C
\end{equation}
where $\bot$ refers to the orthogonal complement.

As example applications, $\mathbf X$ can be a manifold in a functional space or matrix in a vector space, and $\mathbf Y$ a sub manifold or sub matrix of $\mathbf X$.

This retrieval product is clearly not commutative, i.e. the $\mathbf X$-structure in $\mathbf Y$ is not necessarily the same as the $\mathbf Y$-structure in $\mathbf X$. The exception occurs only when $\mathbf X \equiv \mathbf Y$.

\subsection{Decomposition Operator}

Let $\mathbf X$ be a generic $n$-dimensional functional or rank-$n$ tensor spanning an $n$-dimensional space $\mathbb S_\mathbf X$, and let $\mathbf X_\alpha \in \mathbf X$ be generic $k$-dimensional subspace constituents of $\mathbf X$, with every $\alpha \in \{i,\cdots,j\}$, $1 \leqslant i \leqslant j \leqslant n$ and every $k \leqslant n$.

The subspace constituents $\mathbf X_\alpha$ are not necessarily independent from each other. Instead, they can interact through binding functionals $\mathbf C$ of dimension $c \leqslant \min(K)$, where $K$ is the set of $k$ values quantifying the dimension of the subspace constituents taken into consideration.

We hereby define an \textit{Interacting Subspace Decomposition Operator} $\Upsilon_\alpha$(.) as:
\begin{equation}
\Upsilon_\alpha(\mathbf X) \equiv 
\mathbf X \star \mathbf e_\alpha \equiv \mathbf X_\alpha
\end{equation}
where $\mathbf e_\alpha$ is the functional basis of the subspace spanned by $\mathbf X_\alpha$, and $\star$ is the retrieval product defined in equation \eqref{star}.

The application of $\Upsilon(.)$ to a full partition of $\mathbf X$ denoted as $\mathbb P(\mathbf X)$ is then given as the set of all $\mathbf X_\alpha$ within $\mathbb P(\mathbf X)$.
The partitions are generic and not necessarily disjoint.

\section{Interacting Subspace Composition}\label{AC}

\subsection{Composition Product}

Consider $m$ processes represented by generic $k_i$-dimensional functionals $\mathbf X_i$, $i \in \mathbb M$, with $\mathbb M =\{1,\cdots,m\}$. These processes interact in a $c$-dimensional functional $\mathbf C_\mathbb M$, with $c \leqslant \min(k_i)$.
The \textit{Composition Product} $\#$ among the functionals $\mathbf X_i$ is defined as:
\begin{equation}\label{cp}
\mathbf X_1 \, \# \, \cdots \, \# \, \mathbf X_m \equiv \mathbf X_1 \otimes \cdots \otimes \mathbf X_m \odot \mathbf C_\mathbb M
\end{equation}
with
\begin{equation}
\mathbf C_\mathbb M \equiv \mathcal \langle \mathbf X_i,\mathbf X_j \rangle, i,j \in \mathbb M.
\end{equation}

Let $d_{cp}$ denote the dimensionality or rank of the \textit{lhs} (left-hand-side) of \eqref{cp}. Then, $d_{cp} = \sum_p k_p - c$, i.e. the composition product is rank additive, albeit discounting the dimensionality of the interactions.

\subsection{Composition Operator}

Consider now a partition $\mathbb P(\mathbf X)$ into interacting sub manifolds $\mathbf X_\alpha$ as defined in the previous Appendix (section \ref{AB}).
Our aim is to generate the full manifold $\mathbf X$ from its subspace constituents. This can be done by introducing a \textit{Interacting Subspace Composition Operator} (or \textit{Composition Operator} for short) $\Omega(.)$ that essentially implements the Composition Product among the partition members (sub manifolds) of $\mathbf X$:

\begin{equation}
\Omega[\mathbb P(\mathbf X)] \equiv {\mathbf X_\alpha}_{,1} \, \# \, \cdots \, \# {\mathbf X_\alpha}_{,k} = \mathbf X.
\end{equation}
with $\alpha_{,i}$ denoting the $i^{th}$ partition member and $k$ denoting the cardinality of the partition (i.e. the number of its members).

The interacting subspace composition of a spatiotemporal structure generalises the traditional invariant space-time manifolds in Mathematical Physics by taking into account the spatiotemporal codependence quantified by the dynamic codependence functional $\mathbf C$ corresponding to the spatiotemporal coevolution manifold. 
That is, we operate in a generalised setting where space and time are not necessarily independent dimensions in a structurally invariant manifold, but actually coevolving via nonlinear dynamical interactions \citep{Perdigao2016}.

\vspace{12pt}
\paragraph{Acknowledgements:}

\textit{This research was supported by the ERC Advanced Grant 'Flood Change' project no. 291152. C.P. also acknowledges support from the Portuguese Foundation for Science and Technology (FCT) through the project RECI/GEO-MET/0380/2012 - SHARE - Seamless High-resolution Atmosphere-ocean Research, and through the FCT funding contract UID/GEO /50019/2013 for Instituto Dom Luiz. Our families are gratefully acknowledged for their omnipresent support.}

\break

\end{document}